\documentclass{article}%
\usepackage{amsmath}
\usepackage{amsfonts}
\usepackage{amssymb,float}
\usepackage{latexsym}
\usepackage[latin1]{inputenc}
\usepackage{amsmath}
\usepackage{amssymb,float}
\usepackage{latexsym}
\usepackage{epstopdf}
\usepackage{geometry}
\usepackage[active]{srcltx}
\usepackage{graphicx}
\usepackage{epstopdf}
\usepackage{enumerate}
\usepackage[ bookmarks=true,         bookmarksnumbered=true, colorlinks=true,
pdfstartview=FitV, linkcolor=blue, citecolor=blue, urlcolor=blue]{hyperref}
\usepackage{amssymb}%
\providecommand{\U}[1]{\protect\rule{.1in}{.1in}}
\def\R{{\mathbb R}}

\def\E{{{\mathbb E}\,}}
\def\P{{\mathbb P}}

\def\Z{{\mathbb Z}}

\newtheorem{theorem}{Theorem}[section]

\newtheorem{corollary}{Corollary}[section]

\newtheorem{example}{Example}[section]
\newtheorem{lemma}{Lemma}[section]

\newtheorem{remark}{Remark}[section]

\newenvironment{proof}[0]{\paragraph{Proof.}}{\rule{0.5em}{0.5em}}

\begin{document}

\title{ Limit theorems for linear random fields with innovations in the domain of
attraction of a stable law}
\maketitle

\begin{center}
\bigskip

Magda Peligrad$^{a}$, Hailin Sang$^{b}$, Yimin Xiao$^{c}$ and Guangyu
Yang$^{d}$

\bigskip

$^{a}$ Department of Mathematical Sciences, University of Cincinnati,
Cincinnati, OH 45221, USA.\newline E-mail address: peligrm@ucmail.uc.edu

$^{b}$ Department of Mathematics, University of Mississippi, University, MS
38677, USA. \newline  E-mail address: sang@olemiss.edu

\bigskip$^{c}$ Department of Statistics and Probability, Michigan State
University, East Lansing, MI 48824, USA. \newline E-mail address: xiao@stt.msu.edu

\bigskip$^{d}$ School of Mathematics and Statistics, Zhengzhou University,
450001, Henan, China. \newline E-mail address: guangyu@zzu.edu.cn
\end{center}

\bigskip\textbf{Abbreviated Title: }{\Large Limit theorems for linear random
fields}

\begin{center}
\bigskip

\textbf{Abstract}

\end{center}

In this paper we study the convergence in distribution and the local limit
theorem for the partial sums of linear random fields with i.i.d. innovations
that have infinite second moment and belong to the domain of
attraction of a stable law with index $0<\alpha\leq2$ under the condition 
that the innovations are centered if $1<\alpha\leq2$ and are symmetric if
$\alpha=1$. We establish these two types of limit theorems as long as the
linear random fields are well-defined, the coefficients are either absolutely
summable or not absolutely summable.

\noindent\textbf{Keywords:} \noindent Linear random fields,  domain of attraction of a stable law, 
local limit theorem, weak convergence

\bigskip


\section{Introduction}

For a sequence of random variables $\{\xi_{i}\}_{i\in{\mathbb{Z}}^{d}}$ indexed by vectors in $\mathbb{Z}^{d}$,
many statistical procedures produce estimators of the type
\begin{equation}
S_{n}=\sum\limits_{i\in{\mathbb{Z}}^{d}}b_{ni}\xi_{i}, \label{rfsn}%
\end{equation}
where $\{b_{ni}\}_{i\in \Z^d}$ are real numbers. 
To give an example, let us consider the nonlinear regression model
\[
Y(x)=g(x)+\eta(x),
\]
where $g(x)$ is an unknown function and $\eta(x)$ is the random noise. If we fix
the design points $\left\{x_{i}\right\}_{i\in\Gamma_n^d}$, where $\Gamma_n^d$ is a sequence of 
finite regions of $\mathbb{Z}^{d}$, we get
\[
Y_{i}=Y(x_{i})=g(x_{i})+\eta(x_{i})=g(x_{i})+\eta_{i}.
\]
A nonparametric kernel estimator of $g(x)$ is defined as
\begin{equation}\label{Eq:K}
\hat{g}_{n}(x)=\sum\limits_{i\in\Gamma_n^d}\omega_{n,i}(x)Y_{i}
\end{equation}
with
\[
\omega_{n,i}(x)=\frac{K((x_{i}-x)/h_{n})}{\sum\nolimits_{j\in\Gamma_n^d}K((x_{j}%
	-x)/h_{n})},
\]
where $K$ is a kernel function and $h_n$ is a sequence of bandwidths which goes to zero as $n
\rightarrow \infty$.  Define ${g}_{n}(x)=\sum_{i\in\Gamma_n^d}\omega_{n,i}(x)g(x_{i})$; then it is 
clear that, $\hat{g}_{n}(x) - g_n(x)$ is of the type (\ref{rfsn}). Moreover, if the noise 
$\{\eta_{i}\}_{i\in\mathbb{Z}^d}$ is given by a linear random field $\eta_i = \sum_{j\in {\mathbb{Z}}^{d}} 
a_{j}\xi_{i-j}$,  under the conditions that the innovations are independent and identically 
distributed (i.i.d.) random variables with a bounded second moment, the estimator (\ref{Eq:K})  
has been well-studied in the literature, see, e.g., Tran (1990), Hallin, Lu and Tran (2004), El 
Machkoui (2007), El Machkouri and Stoica (2010), and Beknazaryan, Sang and Xiao (2019). 

As we shall see later, another important example of statistics of type
(\ref{rfsn}) appears in the study of the partial sums%
\[
S_{n}=\sum\limits_{i\in\Gamma_n^d}X_{i}
\]
of the versatile linear random field $\{X_{i}\}_{i\in{\mathbb{Z}}^{d}}$, 
expressed as a linear combination
\begin{equation}
	X_{j}=\sum_{i\in{\mathbb{Z}}^{d}}a_{i}\xi_{j-i} \label{deflin}
\end{equation}
with real coefficients $\{a_{i}\}_{i\in\mathbb{Z}^d}$, where $\Gamma_n^d$ is a 
sequence of finite regions of $\mathbb{Z}^{d}$.

Under the condition ${{\mathbb{E}}\,}\xi_{0}^{2}<\infty$, the linear random
fields in the form of (\ref{deflin}) have been intensively studied in the
literature. Among the last decade contributions, we would like to mention
first that Mallik and Woodroofe (2011) studied the CLT for the partial sums of
linear random fields over rectangular regions. Lahiri and Robinson (2016)
established a CLT for the sums over dilated regions. Sang and Xiao (2018)
extended the exact moderate and large deviation limit theorems studied in
Peligrad et al. (2014) from linear processes to linear random fields.
Beknazaryan, Sang and Xiao (2019) established the Cramér type moderate
deviation for the partial sums of linear random fields. Recently, Fortune,
Peligrad and Sang (2021) proved a local limit theorem under some regularity
conditions for the innovations and the sampling regions. The authors of the
above research considered both short range and long range dependence cases for
linear random fields with ${{\mathbb{E}}\,}\xi_{0}^{2}<\infty$. Koul, Mimoto
and Surgailis (2016) studied the goodness-of-fit test for marginal
distribution of linear random fields with long range dependance. For additional work
in studying the limit theorems of linear random fields under the condition
that the innovation has finite second moment, see, for example, the review in
Lahiri and Robinson (2016) and Sang and Xiao (2019) and the references therein.

On the other hand, a rather small number of papers treat the asymptotics of
linear random fields when the innovations do not have finite second moment, in
particular, when the innovations belong to the domain of attraction of a
stable law with index $0<\alpha\leq2$. In the one-dimensional case, regarding the model (\ref{rfsn}),
Shukri (1976) studied the stable limit theorem under some specified conditions
on the weights $\{b_{ni}\}_{i\in \Z}$. Astrauskas (1983) studied the limit
theorem and the functional limit theorem for linear processes in
(\ref{deflin})\ in the case where the coefficients have some specific form.
Davis and Resnick (1985) established a limit theorem for the partial sums of
linear processes with $0<\alpha<2$ when the coefficients are at least
absolutely summable. Balan, Jakubowski and Louhichi (2016) studied the
functional limit theorem for linear processes when the coefficients are
absolutely summable. McElroy and Politis (2003) studied the limit theorem for
the partial sums of linear random fields over one rectangle under the
condition $1<\alpha<2$, the coefficients are summable and the width of each
side of the rectangle goes to infinity.

Our paper will deal with two types of results: convergence in distribution and
local limit theorems for the models (\ref{rfsn})\ and (\ref{deflin}) with
innovation in the domain of attraction of a stable law with index
$0<\alpha\leq2$. A local limit theorem is a more delicate limit result than
the corresponding convergence in distribution result.
More specifically, based on a convergence in distribution theorem in Shukri
(1976), for the model (\ref{rfsn}) we obtain a general local limit theorem
that is applicable to the situation when the innovation's characteristic 
function might not be integrable.

The model (\ref{deflin}) is very difficult to study and, as far as we know,
concerning the local limit theorem, there is no completed work in the
literature even for linear processes (i.e. $d=1$ in our setting).  As a matter of fact,
as mentioned above, the convergence in distribution has not been studied in
its full generality and only partial results have been developed for the
convergence in distribution. 
For the model (\ref{deflin}) we
provide both convergence in distribution and local limit theorem for random
fields. In the case that the coefficients are not absolutely summable, the domain 
of summation is a union of rectangles and our condition involves the rate of 
convergence to infinite of the number of rectangles. Otherwise we only require 
that the normalizer goes to infinity.

We refer to Petrov (1975), Dolgopyat (2016) and the references therein for a 
review on local limit theorem of partial sums of independent random variables.
It is remarkable that, in a recent paper, Fortune, Peligrad and Sang (2021) 
studied the local limit theorem for linear
random fields in the form of (\ref{deflin}) under the conditions that the
innovations have mean zero and finite second moment.

The paper is organized as follows. We first introduce some preliminaries and 
notations in Section 2. To prepare for the local limit theorems we present 
convergence in distribution results for random fields in Section 3. Next, 
Section 4 contains our local limit theorems. We discuss some examples 
in Section 5. Section 6 is dedicated to the proofs.

\section{Notations and preliminaries}
Throughout the paper we shall use the following notations. We use $\iota$ to denote 
$\sqrt{-1}$, the imaginary unit. By $\varphi_{X}(v)$, we denote the characteristic function of
the random variable $X$, i.e.,
\[
\varphi_{X}(v)={{\mathbb{E}}\,}(\exp(\iota vX)).
\]
The notation $\Rightarrow$ is used for convergence in distribution and also
for convergence in probability to a constant. For sequences of positive 
constants $a_n$ and $b_n$, the notation $a_{n}\ll b_{n}$ means that 
$a_{n}/b_{n}$ is bounded; by $a_n \asymp b_n$ we mean that the ratios 
$a_n/b_n$ are bounded from below and above by positive and finite constants. 
$f(r) \sim g(r)$ means that $\lim f(r)/g(r)=1$ as $r\rightarrow \infty$ or $r\rightarrow 0$ 
depending on the context. We denote the cardinality of a  set $\Gamma$ by 
$|\Gamma|$. $c_\alpha$, $C$ and $C_i, 1\le i\le 4$ denote positive constants, which 
are independent of $x$, $t$ or $n$ and may change values from line to line. The 
constant $c_\alpha$ may depend on $\alpha$.

Let $\{\xi_{i}\}_{i\in{\mathbb{Z}}^{d}}$ be i.i.d. random variables in the
domain of attraction of a stable law with index $0<\alpha\leq2$ and
${{\mathbb{E}}\,}\xi_{0}^{2}=\infty$. They satisfy the following condition
\begin{align}\label{tp}
\lim_{x\rightarrow\infty}\frac{\P(|\xi_0|>x)}{x^{-\alpha }L_1(x)}=1,
\end{align}
where $L_1(x)$ is a slowly varying function at infinity.  When $0<\alpha<2$,
\begin{equation}
\frac{{\mathbb{P}}(\xi_{0}>x)}{{\mathbb{P}}(|\xi_{0}|>x)}\rightarrow
c^{+}\;\;\;\text{and}\;\;\;\frac{{\mathbb{P}}(\xi_{0}<-x)}{{\mathbb{P}}%
(|\xi_{0}|>x)}\rightarrow c^{-} \ \ \ \  \text{ as }x\rightarrow\infty. \label{lrtp}%
\end{equation}
Here $0\leq c^{+}\leq1$ and $c^{+}+c^{-}=1$.

By (8.5) of Feller (1971), page 313,  for $0<\alpha\le 2$, 
\begin{align}\label{F83}
\lim_{x\rightarrow\infty}\frac{x^{2}\P(|\xi_0 |>x)}{\E(\xi_0 ^{2}I(|\xi_0 |\leq x))}= \frac{
2-\alpha }{\alpha }.
\end{align}
Combining the two relations (\ref{tp}) and (\ref{F83}) we obtain%
\[
\lim_{x\rightarrow\infty}\frac{x^{2-\alpha }L_1(x)}{\E(\xi_0 ^{2}I(|\xi_0 |\leq x))}=\frac{%
2-\alpha }{\alpha }.
\]
Let  
\begin{equation*}
b=\inf \left\{ x\ge 1: \E(\xi_0 ^{2}I(|\xi_0 |\leq x))>1\right\} .
\end{equation*}%
For $0<\alpha <2$, we define 
\begin{align} \label{defL}
\begin{split}
L(x) &=\frac{2-\alpha }{\alpha }x^{\alpha -2}\E(\xi_0 ^{2}I(|\xi_0 |\leq x))
\;\;\text{ for }x\geq b,\\
\text{ \ and} \;\;\; \;
L(x) &=\frac{2-\alpha }{\alpha }b^{\alpha -2}\E(\xi_0 ^{2}I(|\xi_0 |\leq b))
\;\;\text{ for }0<x<b\text{.}  
\end{split}
\end{align}
Also for $\alpha =2$, we define 
\begin{align} \label{defL2}
\begin{split}
L(x) &=\E(\xi_0 ^{2}I(|\xi_0 |\leq x))\;\;\text{ for }x\geq b,\\
\text{and} 
\;\;\;\;L(x) &=\E(\xi_0 ^{2}I(|\xi_0 |\leq b))\;\;\text{ for }\text0<x<b.
\end{split}
\end{align}
\begin{remark}\label{sv}
By this definition, for all $0<\alpha\le 2$, the function $L(x)$ is slowly varying at infinity,  the function $x^{2-\alpha }L(x)$ is non-decreasing, continuous from the right and has left-hand limits, and $\lim_{x\rightarrow \infty}x^{2-\alpha }L(x)=\infty$. In addition, $\lim_{x\rightarrow\infty}L_1(x)/L(x)=1$, 
$L(x)\geq (\max\{b,x\})^{\alpha -2}(2-\alpha )/\alpha $ for $0<\alpha <2$ and in the case $\alpha =2$, $L(x)\geq 1
$,  $\lim_{x\rightarrow\infty}L_1(x)/L(x)=0$ and $L(x)=(1+o(1))\int_{0}^{x}t^{-1}L_1(t)dt$ 
as $x\rightarrow\infty$. See page 313 of Feller (1971) and page 83 of Ibragimov and Linnik (1971). 
\end{remark}


For the main results of this paper, we assume the following
condition:

\noindent\textbf{Condition }${\boldsymbol A}:$ ${{\mathbb{E}}\,}\xi_{0}=0$ if $1<\alpha
\leq2$ and the innovation $\xi_{0}$ has symmetric distribution if $\alpha=1$.

We apply the following convention throughout the paper.

\noindent\textbf{Convention:} For $x=0$, $p>0$, $|x|^pL(1/|x|)=0$, because $L(y)$ is slowly varying at infinity. 

Under Condition ${\boldsymbol A}$, the characteristic function $\varphi_{\xi}(t)$ of the innovations has
the form 
\begin{align}\label{cha1}
\varphi_{\xi}(t)=\exp \big\{-c_\alpha|t|^\alpha {L}(1/|t|) \big(1 - \iota \beta \tau(\alpha, t)\big)\big\}
\end{align}
for $t$ in the neighborhood of zero, where $c_\alpha>0$, 
\[
\tau(\alpha, t) = \left\{\begin{array}{cc}
{\rm sgn} t \tan\big(\pi \alpha/2\big), \ \ &\hbox{ if }\, \alpha \ne 1,\\
0, &\hbox{ if }\, \alpha = 1,
\end{array}
\right.
\]
and $\beta =c^+ - c^-$ if $0<\alpha<2$. Notice that $\tau(\alpha, t) =0$ if $\alpha=2$.  See, for example, 
Ibragimov and Linnik (1971, Theorem 2.6.5) for the case $\alpha\neq1$, Aaronson and Denker
(1998, Theorem 2) for $\alpha=1$. 

We shall introduce some conditions on the innovations for the local limit
theorem. Recall that a random variable $X$ does not have a lattice distribution if and only if
$|\varphi_{X}(t)|<1$ for all $t\neq0$. On the other hand, the Cramér
condition means that $\lim\sup_{|t|\rightarrow\infty}|\varphi_{X}(t)|<1$. It should be mentioned 
that $X$ has a non-lattice distribution whenever it satisfies the Cramér condition. 
The Cramér condition was also called strongly non-lattice in Stone (1965).

For comparison purpose, let us recall some known results for linear processes. For i.i.d. random 
variables $\{\xi_{i}\}_{i=1}^{\infty}$, which satisfy
conditions (\ref{tp}), (\ref{lrtp}) and Condition $\boldsymbol A$, let $\sigma_{n}$ be a sequence of
positive numbers such that
\[
n{\mathbb{P}}(|\xi_{1}|>\sigma_{n}x)\rightarrow x^{-\alpha}\;\;\text{as}%
\;\;n\rightarrow\infty\;\;\text{for all}\;\;\;x>0\text{ and }0<\alpha<2.
\]
The normalizing constant $\sigma_{n}$ can also be defined as $\inf\{x:{\mathbb{P}}(|\xi_{1}|>x)\leq
n^{-1}\}$. 
In general, up to a constant factor ($(2-\alpha)/{\alpha}$ if $0<\alpha
<2$, see Feller, 1971), $\sigma_{n}^{2}$ can be chosen to satisfy%
\[
\frac{n}{\sigma_{n}^{2}}{{\mathbb{E}}\,}\xi_{1}^{2}I(|\xi_{1}|\leq\sigma
_{n})\rightarrow1,\text{ as }n\rightarrow\infty,0<\alpha\leq2.
\]
Then, by  Feller (1971, p. 580), we have the following limit theorem
\begin{equation}\label{Eq:stablelimt}
\sigma_{n}^{-1}\sum_{i=1}^{n} \xi_{i} \, \Rightarrow \,S
\end{equation}
as $n\rightarrow\infty$, where $S$ is an $\alpha$-stable random variable  
with characteristic function 
\begin{align}\label{cs}
{{\mathbb{E}}\,}(\exp(\iota t S)) = \exp\big\{  -c_\alpha |t|^\alpha \big(1 - \iota \beta \tau(\alpha, t)\big)\big\}
\end{align}
with the same $c_\alpha$, $\beta$ and  $\tau(\alpha, t)$ as in \eqref{cha1}. In the sequel, we will denote the distribution function of $S$ by $G(\cdot)$.

Under the conditions (\ref{tp}) and (\ref{lrtp}), for one-sided linear
processes $X_{j}=\sum_{i=0}^{\infty}a_{i}\xi_{j-i}$, $S_{n}=\sum_{j=1}%
^{n}X_{j}$, Davis and Resnick (1985) assumed that the real coefficients
satisfy
\[
\sum\nolimits_{i=0}^{\infty}|a_{i}|^{\delta}<\infty\;\;\text{for
some}\;\;\;\delta<\alpha, \;\delta\leq1,
\]
and proved that
\begin{equation*}
\sigma_{n}^{-1}\left(S_{n}-nC_{n}\sum_{i=0}^{\infty}a_{i}\right)\Rightarrow\left(\sum
_{i=0}^{\infty}a_{i}\right)\tilde{S} \label{DR}%
\end{equation*}
as $n\rightarrow\infty$, where $C_{n}=0$ if $0 < \alpha \le 1$ and $C_{n} = \mathbb{E}\big[\xi_{1}I(|\xi_{1}|\leq\sigma_{n})\big]$
if $1 <\alpha <2$, and where $\tilde{S}$ is a stable random variable with index
$\alpha$.  If $\sum_{i=0}^{\infty}a_{i}=0$, then $\sigma_{n}%
^{-1}S_{n}\Rightarrow0$. We remark that if, in addition, Condition $\boldsymbol A$ holds, then $\tilde{S}$ is the 
random variable $S$ in (\ref{Eq:stablelimt}).


\section{Convergence in distribution for random fields}

We mention that we can easily extend the results in Shukri (1976) to triangular arrays of random fields.

Let $\{c_{ni}\}_{i\in{\mathbb{Z}}^{d}}$ be real numbers. We shall introduce two
conditions:\bigskip

\textbf{Condition} ($A_{1}$): $\sum\limits_{i\in{\mathbb{Z}}^{d}}|c_{ni}%
|^{\alpha}L(1/|c_{ni}|)\rightarrow1$ as $n\rightarrow\infty$, where $L(\cdot)$ is defined in (\ref{defL}) or (\ref{defL2}) \newline and

\textbf{Condition} ($A_{2}$): $\rho_{n}=\sup_{i}|c_{ni}|\rightarrow0$ as
$n\rightarrow\infty$. \bigskip

We consider a sequence  $\{S_{n}\} $ of random variables with the form 
$S_{n}=\sum_{i\in{\mathbb{Z}}^{d}}c_{ni}\xi_{i}$. This model is a specialized
case of model (\ref{rfsn}) if we impose additional conditions on the
coefficients $\{c_{ni}\}$ such as Condition ($A_{1}$) or ($A_{2}$). With a similar
proof as in Theorem 1, Corollary 1 and Corollary 2 in Shukri (1976), we
formulate the following general result and mention one corollary. 

\begin{theorem}
\label{series} 
For $0<\alpha\leq2$, assume that innovations $\{\xi_{i}\}_{i\in{\mathbb{Z}}^{d}}$ 
are i.i.d. random variables in the domain of attraction of an $\alpha$-stable law.  The innovations satisfy Condition ${\boldsymbol A}$ and have characteristic function (\ref{cha1}).  The coefficients $\{c_{ni}\}$ satisfy Condition ($A_{2}$) 
and
\begin{equation}\label{Eq:S1}
\lim_{n\rightarrow \infty}\sum\limits_{i\in{\mathbb{Z}}^{d},c_{ni}>0}|c_{ni}|^{\alpha}L(1/|c_{ni}|%
)=c
\end{equation}
\begin{equation}\label{Eq:S2}
\lim_{n\rightarrow \infty}\sum\limits_{i\in{\mathbb{Z}}^{d},c_{ni}<0}|c_{ni}|^{\alpha}L(1/|c_{ni}|%
)=1-c,
\end{equation}
for some $0\leq c\leq1$. Then
the sequence $S_{n} = \sum_{i\in{\mathbb{Z}}^{d}}c_{ni}\xi_{i}$ converges
weakly to an $\alpha$-stable random variable with the form $ c^{1/\alpha}S' - (1-c)^{1/\alpha} S''$, where $S'$ and $S''$ are independent $\alpha$-stable 
random variables that have distribution function $G(\cdot)$ and characteristic function (\ref{cs}).
\end{theorem}

Notice that the probability distribution function of 
the limiting random variable in Theorem \ref{series}  can be represented as 
$G\big( c^{-1/\alpha}\, \cdot \big)\ast \big(1-G(- (1-c)^{-1/\alpha}\, \cdot)\big)$ for $0<c<1$; it is $G(x)$ for $c=1$ and $1-G(-x)$ for $c=0$. 
The following 
is a direct consequence of Theorem \ref{series}.

\begin{corollary}
\label{series1} 
For $0<\alpha\leq2$, assume that innovations $\{\xi_{i}\}_{i\in{\mathbb{Z}}^{d}}$ 
are i.i.d. random variables in the domain of attraction of an $\alpha$-stable law.  The innovations satisfy Condition ${\boldsymbol A}$ and have characteristic function (\ref{cha1}).  The coefficients $\{c_{ni}\}$ satisfy Conditions ($A_{1}$) and ($A_{2}$). 
If either the innovations have symmetric distribution or $c_{ni}\geq0$ for all $n$ and $i$, then 
$S_{n}$ converges weakly to a random variable that has distribution function $G(\cdot)$ and characteristic function (\ref{cs}).
\end{corollary}



Next, for the sake of applications we shall first study the normalized form of
these results and provide more flexible sufficient conditions on the real
coefficients. These results will be especially useful for analyzing the
partial sums of linear random fields $\{X_{j}\}_{j\in{\mathbb{Z}}^{d}}$ defined
by (\ref{deflin}), namely%
\[
X_{j}=\sum_{i\in{\mathbb{Z}}^{d}}a_{i}\xi_{j-i}.
\]
By the three-series theorem, the linear random field in (\ref{deflin})
converges almost surely if and only if
\begin{align}\label{existence}
\sum_{i\in{\mathbb{Z}}^{d}}|a_{i}|^{\alpha}L(1/|a_{i}|)<\infty.
\end{align}
For the one-dimensional case $d=1$, this statement can be found in Astrauskas (1983) and in Proposition 5.4 in Balan, Jakubowski and Louhichi (2016) for the case
$0<\alpha<2$ and in the proof of Theorem 2.1 in Peligrad and Sang (2012) for the case $\alpha=2$. The proof for the general $d$-dimensional linear random
fields is similar. With the same argument, 
the linear random field of the form (\ref{rfsn}) converges almost surely if and only if for any $n$
\begin{align}\label{existence1}
\sum_{i\in{\mathbb{Z}}^{d}}|b_{ni}|^{\alpha}L(1/|b_{ni}|)<\infty.
\end{align}
For all the results in this paper, we assume condition (\ref{existence}) for linear random field in (\ref{deflin}) and  condition (\ref{existence1}) for linear random field in (\ref{rfsn}). 

Let us normalize the linear random fields in (\ref{rfsn}). Recall the definition of $L(x)$ in (\ref{defL}) and (\ref{defL2}). 
For $0<\alpha\leq 2$ define
\begin{align}\label{defBn}
B_{n}=\inf\left\{  x\ge 1:\sum_{i\in\mathbb{Z}^d}(|b_{ni}|/x)^{\alpha}L(x/|b_{ni}|)\leq
1\right\}  \text{.}
\end{align}
By the definition of $L(x)$, it can be shown that $B_{n}$ is a sequence of positive numbers such that
\begin{equation}
\sum_{i\in\mathbb{Z}^d}(|b_{ni}|/B_{n})^{\alpha}L(B_n/|b_{ni}|)=1.\label{limit}
\end{equation}
See the proof in Lemma \ref{bn}. In the case $\alpha=2,$ and ${{\mathbb{E}%
}\,}\xi_{0}^{2}=\infty$, Peligrad and Sang (2013) used the same definition
for the normalizer $B_{n}$ for studying a triangular array of martingale
differences.

For each $0<\alpha\leq2$, with $c_{ni}%
=b_{ni}/B_{n}$, Condition ($A_{1}$) always holds because of (\ref{limit}). Under the assumption 
$B_{n}\rightarrow\infty$ and $\{b_{ni}\}$ are bounded uniformly on $n$ and $i$, Condition 
($A_{2}$) holds for
$c_{ni}=b_{ni}/B_{n}$ and therefore Theorem \ref{series} and Corollary
\ref{series1} 
hold for linear random fields (\ref{rfsn}).

If $\limsup_{n}\sup_{i\in\mathbb{Z}^{d}}|b_{ni}|=\infty$, in Lemma
\ref{lemmarho}, we shall provide sufficient conditions for Condition ($A_{2}$),
which are easier to verify. We first introduce the following notations, as in
Fortune, Peligrad and Sang (2021); see also Mallik and Woodroofe (2011). For a
countable collection of real numbers $\{b_{i}, i\in{\mathbb{Z}}^{d}\}$, where
$i=(i_{1},...,i_{d})$, we denote an increment in the direction $k$ by
\[
\Delta_{k}b_{i_{1},...,i_{k},...,i_{d}}=b_{i_{1},...,i_{k},...,i_{d}}%
-b_{i_{1},...,i_{k}-1,...,i_{d}}%
\]
and their composition is denoted by $\Delta:$
\begin{equation}
\Delta b_{i}=\Delta_{1}\circ\Delta_{2}\circ...\circ\Delta_{d}b_{i}.
\label{Incr}%
\end{equation}
For a collection of sequences $\{b_{ni}, i\in{\mathbb{Z}}^{d}\}$, define $\Delta b_{ni}$ 
by (\ref{Incr}) with $b_i$ replaced by $b_{ni}$ and denote
\[
\Delta_{n}=\sup_{i\in\Z^d}|\Delta b_{ni}|.
\]
Define
 \[
\rho_{n}=\sup_{i\in\Z^d}|b_{ni}|/B_{n}.
\]
By the first inequality in Lemma \ref{lemmarho}, $\rho_{n}\rightarrow0$ as
$n\rightarrow\infty$ if $\Delta_{n}=o(B_{n})$. Together with the above results
and analysis, we have the following limit theorem for the linear model (\ref{rfsn}).

\begin{theorem}\label{thm1'} 
For the linear model (\ref{rfsn}) with i.i.d. innovations $\{\xi_{i}\}_{i\in{\mathbb{Z}}^{d}}$ in the domain of attraction of 
an $\alpha$-stable law ($0<\alpha\le 2$), assume that the innovations satisfy Condition ${\boldsymbol A}$ and have characteristic function (\ref{cha1}), the $B_{n}$ defined in
(\ref{defBn}) satisfies $B_{n}\rightarrow\infty$ as
$n\rightarrow\infty$. In the case that $\limsup_{n\to \infty}\sup_{i\in\mathbb{Z}^{d}
}|b_{ni}|=\infty$, we additionally require that $\Delta_{n}%
=o(B_{n})$.  
\begin{enumerate}
\item If conditions (\ref{Eq:S1}) and (\ref{Eq:S2}) hold for $c_{ni} =b_{ni}/B_{n}$, then $S_{n}/B_n$ converges
weakly to an $\alpha$-stable random variable with the form $ c^{1/\alpha}S' - (1-c)^{1/\alpha} S''$, where $S'$ and $S''$ are independent $\alpha$-stable 
random variables that have distribution function $G(\cdot)$ and characteristic function (\ref{cs}).
\item If either  the innovations have symmetric distribution or the coefficients $b_{ni} \ge 0$
for all $n$ and $i$, then $S_{n}/B_n$ converges
weakly to an $\alpha$-stable random variable with distribution function $G(\cdot)$. 
\end{enumerate}
\end{theorem}

Now we study the linear random fields (\ref{deflin}).
Let $\Gamma_{n}^{d}$ be a sequence of finite subsets of ${\mathbb{Z}}^{d}$,
and define the sum
\begin{equation}
S_{n}=\sum_{j\in\Gamma_{n}^{d}}X_{j}, \label{S_n}%
\end{equation}
where $X_{j}$ is a linear random field given in (\ref{deflin}). With the
notation
\begin{equation}
b_{ni}=\sum\limits_{j\in\Gamma_{n}^{d}}a_{j-i}, \label{formb}%
\end{equation}
the sum $S_{n}$ can be expressed as an infinite linear combination of the
innovations as in form (\ref{rfsn}). We still have the definition of $B_n$ as in (\ref{defBn}). 

Similarly, by the analysis for linear model (\ref{rfsn}) with $c_{ni}%
=b_{ni}/B_{n}$, Condition ($A_{1}$) always holds because of (\ref{limit}). In
the case that $0<\alpha<1$, the existence of (\ref{deflin}) implies that
$\sum_{i\in{\mathbb{Z}}^{d}}|a_{i}|$ is finite and then $\{b_{ni}\}$ are
bounded uniformly on $n$ and $i$. Then Condition ($A_{2}$) holds for
$c_{ni}=b_{ni}/B_{n}$ if we assume $B_{n}\rightarrow\infty$. Therefore with the
condition $B_{n}\rightarrow\infty$, Theorem \ref{series} and Corollary
\ref{series1} 
hold for the partial sums of
linear random fields (\ref{deflin}) in (\ref{S_n}) if $0<\alpha<1$ or $1\leq\alpha\leq2$ and $\sum_{i\in
{\mathbb{Z}}^{d}}|a_{i}|<\infty$. 
We remark that,  in order to have the convergence theorems in these two cases,
except assuming $|\Gamma_{n}^{d}|\rightarrow\infty$ as
$n\rightarrow\infty$ (which is due to
the requirement $B_{n}\rightarrow\infty$), we do not impose any other
restriction on the shape of the regions $\Gamma_{n}^{d}$.

If $1\leq\alpha\leq2$ and $\sum_{i\in {\mathbb{Z}}^{d}}|a_{i}|=\infty$, we require the
index sets to be of the form
\begin{equation}
\Gamma_{n}^{d}=\bigcup\limits_{k=1}^{J_{n}}\Gamma_{n}^{d}(k), \label{defgamma}%
\end{equation}
where $\{\Gamma_{n}^{d}(k)\}_{k=1}^{J_{n}}$ is a pairwise disjoint family of
discrete rectangles in ${\mathbb{Z}}^{d}$ and $\Gamma_{n}^{d}(k)$ has the form
\[
\prod_{\ell=1}^{d}[\underline{n}_{\ell}(k),\overline{n}_{\ell}(k)]\cap
{\mathbb{Z}}^{d}%
\]
with $\underline{n}_{\ell}(k)\leq\overline{n}_{\ell}(k)$, where $1\leq\ell\leq
d$, $1\leq k\leq J_{n}$.

By the second inequality in Lemma \ref{lemmarho},  notice that $J_{n}=o(B_{n}^{q})$ implies $\rho_{n}%
\rightarrow0$ as $n\rightarrow\infty$. Hence if $1\leq\alpha\leq2$ and
$\sum_{i\in{\mathbb{Z}}^{d}}|a_{i}|=\infty$, under the condition
$J_{n}=o(B_{n}^{q})$, the generalization of Theorem 1, Corollary 1 and  Corollary 2 in
Shukri (1976), i.e., Theorem \ref{series} and Corollary \ref{series1} hold. 

In summary, we have the following limit theorem for the partial sums of linear random fields (\ref{deflin}).

\begin{theorem}
\label{thm1} 
For the linear random field (\ref{deflin})  with i.i.d. innovations $\{\xi_{i}\}_{i\in{\mathbb{Z}}^{d}}$ in the domain of attraction of 
an $\alpha$-stable law ($0<\alpha\le 2$), let $S_{n}$ be the partial sum defined in 
(\ref{S_n}) and $B_{n}$ be defined as in (\ref{defBn}). Assume that the innovations satisfy Condition ${\boldsymbol A}$ and have characteristic function (\ref{cha1}), and $B_{n}\rightarrow\infty$
as $n\rightarrow\infty$. In the case that $1\le\alpha\le 2$ and $\sum_{i\in{\mathbb{Z}%
}^{d}}|a_{i}|=\infty$, we additionally require that  the
sets $\Gamma_{n}^{d}$ are constructed as a disjoint union of $J_{n}$ discrete
rectangles as in (\ref{defgamma}), where $J_{n}=o(B_{n}^{q})$, $1/p+1/q=1$,
for some $p>\alpha$ if $1\le\alpha<2$, and $p=2$ if $\alpha=2$. Otherwise no
such restriction is required. 
\begin{enumerate}
\item If conditions (\ref{Eq:S1}) and (\ref{Eq:S2}) hold for $c_{ni} =b_{ni}/B_{n}$, where $\{b_{ni}\}$ are defined by \eqref{formb}, then $S_{n}/B_n$ converges
weakly to an $\alpha$-stable random variable with the form $ c^{1/\alpha}S' - (1-c)^{1/\alpha} S''$, where $S'$ and $S''$ are independent $\alpha$-stable 
random variables that have distribution function $G(\cdot)$ and characteristic function (\ref{cs}).
\item If either  the innovations have symmetric distribution or the coefficients $a_{i} \ge 0$, 
for all $i\in\mathbb{Z}^{d}$, then $S_{n}/B_n$ converges
weakly to an $\alpha$-stable random variable with distribution function $G(\cdot)$. 
\end{enumerate}
\end{theorem}
\begin{remark}\label{alpha1}
 In Theorem \ref{thm1},  for  $1\le\alpha<2$, we can take $p>\alpha$ but
arbitrarily close to $\alpha$. Then $q>2$ can be arbitrarily close to $\alpha/(\alpha-1)$. If $\alpha=1$, 
the value of $q$ can be taken arbitrarily large. 
\end{remark}

Let us compare our results with some available results in the literature.

In Theorem 2.1 of McElroy and Politis (2003), the authors only obtained the
convergence theorem for the partial sums of linear random field over one
rectangle under the conditions: $1<\alpha<2$, the coefficients $\{a_{i}\}$ are
summable and $\min_{i}n_{i}\rightarrow\infty$, where $n_{i}$ is the size of
the rectangle in the $i$-th dimension. In the present paper, we
do not have such a restriction on the regions of the sums provided $0<\alpha<1$, or
$1\le\alpha\leq2$ and $\sum\nolimits_{i\in{\mathbb{Z}}^{d}}|a_{i}|<\infty$. The only condition
for the weak convergence is $B_{n}\rightarrow\infty$ as $n\rightarrow\infty$.
If $1\le\alpha\leq2$ and $\sum\nolimits_{i\in{\mathbb{Z}}^{d}}|a_{i}|=\infty$,
we 
obtain the limit theorems for the sums over a finite number $J_{n}$ of
non-empty pairwise mutually exclusive rectangles, as long as $J_{n}%
=o(B_{n}^{q})$ for some $q\geq2$ that is related to $\alpha$. For our
limit results, we do not require the coefficients to be summable if the
existence of the linear random field itself does not implies it, as in the
case $0<\alpha<1$.

The result in Theorem \ref{thm1} is new even in one-dimensional case. Davis
and Resnick (1985) studied the limit theorem for the partial sums $S_{n}%
=\sum_{j=1}^{n}X_{j}$ of one-sided linear process under the condition the
coefficients are summable. In Theorem \ref{thm1}, we do not require the
coefficients to be summable. The sum is not necessarily of the form
$S_{n}=\sum_{j=1}^{n}X_{j}$. If $1\le\alpha\leq2$ and the coefficients are not
summable, we require that $J_{n}=o(B_{n}^{q})$. Otherwise, we only require
that $B_{n}\rightarrow\infty$ as $n\rightarrow\infty$.


\section{Local limit theorems}


\subsection{General local limit theorem}

First we prove a general local limit theorem based on characteristic functions.  
Let $\{S_{n}\}$ be a sequence of random variables 
and $\{B_{n}\}$ be a sequence of positive numbers.

\begin{theorem}\label{general}
Assume that
\begin{align}
\label{weak}\frac{S_{n}}{B_{n}}\Rightarrow \mathcal{L}\text{, }\mathcal{L}\text{ has an integrable
characteristic function and }B_{n}\rightarrow\infty\text{.}%
\end{align}
In addition, suppose that for each $D>0$
\begin{equation}
\ \lim_{T\rightarrow\infty}\text{\ }\lim\sup_{n\rightarrow\infty}%
\int\nolimits_{T<|t|\leq DB_{n}}\left|  {\mathbb{E}}\exp \left(\iota t\frac{S_{n}}{B_{n}%
}\right)\right|  dt=0\text{.} \label{integral}%
\end{equation}
Then for any continuous function $m$ on $\mathbb{R}$ with compact
support,
\begin{equation}
\lim_{n\rightarrow\infty}\sup_{u\in\mathbb{R}}\left|  B_{n}{\mathbb{E}}%
m(S_{n}+u)-f_{\mathcal{L}}\left( -\frac{u}{B_{n}}\right)  \int m(t)\lambda(dt)\right|
=0, \label{conclusion}%
\end{equation}
where $\lambda$ is the Lebesgue measure and $f_{\mathcal{L}}$ is the density of $\mathcal{L}$.
\end{theorem}

\begin{remark}
Since $\mathcal{L}$ has an integrable characteristic function, by Billingsley (1995)
page 371 we know that the limiting distribution has a density.
\end{remark}

\begin{remark}
By decomposing the integral in (\ref{integral}) into two parts, on
$\{T\leq|u|\leq\delta B_{n}\}$ and $\{\delta B_{n}\leq|u|\leq DB_{n}\}$, and
changing the variable in the second integral, we easily argue that in order to
prove this theorem, it is enough to show that for each $D$ fixed there is
$0<\delta<D$ such that%
\[
(D_{1})\text{ \ \ \ }\lim_{T\rightarrow\infty}\lim_{n\rightarrow\infty}%
\sup\int\nolimits_{T<|t|\leq B_{n}\delta}\left|  {\mathbb{E}}\exp
\left(\iota t\frac{S_{n}}{B_{n}}\right)\right|  dt=0
\]
and
\[
(D_{2})\text{ \ \ \ \ \ \ \ }\lim_{n\rightarrow\infty}B_{n}\int
\nolimits_{\delta<|t|\leq D}|{\mathbb{E}}\exp(\iota tS_{n})|dt=0.\text{ \ \ \ \ \ }%
\]

\end{remark}

By taking $u=0$ in (\ref{conclusion}), as in the proof of the classical
Portmanteau weak convergence theorem (Theorem 25.8 in Billingsley, 1995), we obtain:

\begin{corollary}
Under the conditions of Theorem \ref{general} we also have for all $a<b$,
\[
\lim_{n\rightarrow\infty}B_{n}{\mathbb{P}}(a<S_{n}\leq b)=f_{\mathcal{L}}(0)(b-a).
\]

\end{corollary}

\subsection{Local limit theorem for linear random fields}

We first consider the normalized form of the general linear random field
(\ref{rfsn}). This field has the form
$S_{n}=\sum_{i\in{\mathbb{Z}}^{d}}b_{ni}\xi_{i}$ where $\{b_{ni}\}$ is a general sequence of constants 
and does not necessarily has the form (\ref{formb}).
We define $B_{n}$ as in (\ref{defBn}). 

Recall that $\rho_n=\sup_{i\in\mathbb{Z}%
^{d}}|b_{ni}|/B_{n}$. We have the following local limit theorem for model (\ref{rfsn}).
\begin{theorem}
\label{thm2} Assume that the linear random field (\ref{rfsn}) has  i.i.d. innovations $\{\xi_{i}\}_{i\in{\mathbb{Z}}^{d}}$
in the domain of attraction of a stable law with index $0<\alpha\le 2$, 
$B_{n}\rightarrow\infty$ and $\rho_n\rightarrow0$ as $n\rightarrow\infty$. The innovations satisfy 
Condition $\boldsymbol A$ and have a non-lattice distribution. If $\limsup_{n\to \infty}\sup_{i\in\mathbb{Z}^{d}%
}|b_{ni}|=\infty$, we further assume that the innovations satisfy the Cramér
condition. 
Otherwise these additional assumptions are not required. Then, under the conditions for each of the two cases in Theorem \ref{thm1'}, for any continuous function 
$m$ on $\mathbb{R}$ with compact support, 
(\ref{conclusion}) holds with $\mathcal{L}$ being the limiting $\alpha$-stable random variable in Theorem \ref{thm1'}.
\end{theorem}

Notice that if $\limsup\limits_{n\to \infty}\sup\limits_{i\in \mathbb{Z}^{d}}|b_{ni}|<\infty$,  $\rho_n\rightarrow0$ automatically 
since $B_n\rightarrow \infty$. In the case that $\limsup\limits_{n\to \infty}\sup\limits_{i\in \mathbb{Z}^{d}}|b_{ni}|$ $=\infty$,
 by Lemma \ref{lemmarho},  we can replace 
the condition $\rho_n\rightarrow0$ by $\Delta_{n}=o(B_{n})$.

As an application of Theorem \ref{thm2}, we derive the following local limit theorems for the linear
random fields in (\ref{deflin}) with the short range dependence and long range dependence, respectively.

\begin{theorem}
\label{thm3} For the linear random field (\ref{deflin})  with i.i.d. innovations $\{\xi_{i}\}_{i\in{\mathbb{Z}}^{d}}$ in the domain of attraction of 
an $\alpha$-stable law ($0<\alpha\le 2$), let $S_{n}$ be the partial sum defined in 
(\ref{S_n}) and $B_{n}$ be defined as in (\ref{defBn}). Assume that the innovations satisfy Condition $\boldsymbol A$ and have a non-lattice distribution, 
$\sum \nolimits_{i\in{\mathbb{Z}}^{d}}|a_{i}|<\infty$ and  $B_{n}\rightarrow\infty$ 
as $n\rightarrow\infty$.  
Then, under the conditions for each of the two cases in Theorem \ref{thm1}, for any continuous function $m$ on $\mathbb{R}$ with compact support,
(\ref{conclusion}) holds with $\mathcal{L}$ being the limiting $\alpha$-stable random variable in Theorem \ref{thm1}.
\end{theorem}

The next theorem deals with the case when the sum of the coefficients are divergent. In the case where 
$0<\alpha<1$, the existence of the linear random fields implies $\sum \nolimits_{i\in{\mathbb{Z}}^{d}}|a_{i}|
<\infty$, a case which we have already treated it in the above theorem.  Recall that $\xi$ has a non-lattice distribution 
whenever it satisfies the Cram\'er condition.

\begin{theorem}
\label{thm3'} For the linear random field (\ref{deflin})  with i.i.d. innovations $\{\xi_{i}\}_{i\in{\mathbb{Z}}^{d}}$ in the domain of attraction of 
an $\alpha$-stable law ($1\le\alpha\le 2$), let $S_{n}$ be the partial sum defined in 
(\ref{S_n}) and $B_{n}$ be defined as in (\ref{defBn}). Assume that the innovations satisfy Condition $\boldsymbol A$ and the Cramér condition,  $\sum \nolimits_{i\in{\mathbb{Z}}^{d}}|a_{i}|=\infty$ and  $B_{n}\rightarrow\infty$ 
as $n\rightarrow\infty$.  The sum $S_{n}$ is over $J_{n}$
non-empty pairwise mutually exclusive rectangles in $\mathbb{Z}^{d}$ with
$J_{n}=o(B_{n}^{q})$, $1/p+1/q=1$, for some $p>\alpha$ if $1\le \alpha<2$, and $p=2$ of $\alpha=2$. 
Then, under the conditions for each of the two cases in Theorem \ref{thm1},  for any
function $m$ on $\mathbb{R}$, which is continuous and has compact support,
(\ref{conclusion}) holds  with $\mathcal{L}$ being the limiting $\alpha$-stable random variable in Theorem \ref{thm1}.
\end{theorem}


\begin{remark}
The two theorems give us local limit theorems for linear random fields with
innovations in the domain of attraction of a stable law under some regularity
conditions for the innovations and the coefficients as long as the linear
random fields are well defined. In particular, if $1\le\alpha\le2$ and the
coefficients are not absolutely summable, the local limit theorem holds for
the sums over a sequence of regions $\Gamma_{n}^{d}$ which is a disjoint union
of discrete rectangles. In application it allows us to have disjoint discrete
rectangles as spatial sampling regions, and the number of these disjoint
rectangles may increase as the sample size increases. In particular we may
have $(\prod_{k=1}^{d}[\underline{n}_{k},\overline{n}_{k}])\cap{\mathbb{Z}%
}^{d}$ where $\underline{n}_{k}=\overline{n}_{k}$ for some $k$'s. So we may have
a single point region if the equality holds for all $k$'s. The local limit
theorem is new also for $d=1$. Furthermore, we have the freedom to take the
sum over $J_{n}$ blocks of random variables as long as it satisfies the
condition $J_{n}=o(B_{n}^{q})$ for $q\ge2$.
\end{remark}

In fact, with the method of proof as in Theorem \ref{thm2} and Theorem
\ref{thm3}, we can improve the local limit theorem Theorem 2.2 in Fortune,
Peligrad and Sang (2021) for linear random fields when the innovations have
finite variance. Condition (9) there can be improved to $J_{n}=o(B_{n}%
^{2})$. For completeness, we state the local limit theorem in this case below.

\begin{theorem}
\label{thm4} For the linear random field defined in (\ref{deflin}), assume
that the innovations $\xi_{i}$ are i.i.d. random variables with mean zero
$({{\mathbb{E}}\,}\xi_{i}=0)$, finite variance $({{\mathbb{E}}\,}\xi_{i}%
^{2}=\sigma_{\xi}^{2})$, and non-lattice distribution and the collection
$\{a_{i}:i\in{\mathbb{Z}}^{d}\}$ of real coefficients satisfies
\begin{equation*}
\sum\nolimits_{i\in{\mathbb{Z}}^{d}}a_{i}^{2}<\infty. \label{cond coef}%
\end{equation*}
Let $S_{n}$ be defined as in (\ref{S_n}) or its equivalent form (\ref{rfsn})
and assume that $B_{n}=\mathrm{Var}(S_{n})=\sigma_{\xi}^{2}\sum_{i\in
{\mathbb{Z}}^{d}}b_{ni}^{2}\rightarrow\infty$. In the case $\sum
\nolimits_{i\in{\mathbb{Z}}^{d}}|a_{i}|=\infty$, the field has long range
dependence, we additionally assume that the innovations satisfy the Cramér
condition and that the sets $\Gamma_{n}^{d}$ are constructed as a disjoint
union of $J_{n}$ discrete rectangles with $J_{n}=o(B_{n}^{2})$. Under these
conditions, for any function $m$ on $\mathbb{R}$ which is continuous and has compact
support,
\[
\lim_{n\rightarrow\infty}\sup_{u\in\mathbb{R}}\bigg|\sqrt{2\pi}B_{n}%
{{\mathbb{E}}\,}m(S_{n}+u)-[\exp(-u^{2}/2B_{n}^{2})]\int m(x)\lambda
(dx)\bigg|=0,
\]
where $\lambda$ is the Lebesgue measure.
\end{theorem}


\section{Examples}
In this section, we provide some examples to illustrate applications of the main results of this paper.  

\begin{example} {\rm (}{\it Doubly geometric spatial autoregressive models}{\rm)}
{\rm Suppose that $d=2$  and $\{\xi_{i_1,i_2}, (i_1,i_2)\in\mathbb{Z}^2\}$ 
are i.i.d. random variables in the domain of attraction of a stable law with index $0<\alpha\leq2$. 
For any $i=(i_1,i_2)\in\mathbb{Z}^2$, let
\begin{equation}\label{Eq:double}
X_i=X_{i_1,i_2}=\theta X_{i_1-1,i_2}+\rho X_{i_1,i_2-1}-\theta\rho X_{i_1-1,i_2-1}+\xi_{i_1,i_2},
\end{equation}
where $|\theta|<1$ and $|\rho|<1$. By using the back-shift operators $\mathbf{B}_1$ and $\mathbf{B}_2$ 
respectively in the horizontal and vertical directions defined by
\[
\mathbf{B}_1\xi_{i_1,i_2} = \xi_{i_1-1,i_2}, \quad \mathbf{B}_2\xi_{i_1,i_2} = \xi_{i_1,i_2-1}, 
\] 
we have
\begin{align*}
X_{i_1,i_2}&={(1-\theta\mathbf{B}_1)^{-1}(1-\rho\mathbf{B}_2)^{-1}}\xi_{i_1,i_2}\\
&=\sum_{j_1=0}^\infty\sum_{j_2=0}^\infty\theta^{j_1}\rho^{j_2}\xi_{i_1-j_1,i_2-j_2}.
\end{align*}
Thus, $\{X_i, i \in \mathbb Z^2\}$ is a linear random field with  coefficients $a_{j_1,j_2}=\theta^{j_1}\rho^{j_2}$ 
($j_1, j_2\geq 0$) which are summable. When the innovations have finite second moments, this model was introduced by 
Martin (1979) as an analogue in the spatial setting to the classical AR(1) time series. For more information on this and 
more general stationary spatial MA and ARMA models with second moments and their statistical inference one can refer to 
Tj\o{}stheim (1978, 1983),  Guyon (1995), Gaetan and Guyon (2010), Beran et al. (2013). 

Extension of the spatial model (\ref{Eq:double}) to the setting of heavy-tailed innovations can be justified by the 
convergence criterion as in (\ref{existence}).

For any integer $n \ge 1$, let $\Gamma_n^2=[0,n]^2\cap\Z^2$.  
Then for every $i = (i_1, i_2) \in \mathbb Z^2$, 
\begin{align*}
b_{ni}&=\sum_{(j_1,j_2)\in\Gamma_n^2} a_{j_1-i_1,j_2-i_2}
\\&=\left\{
\begin{array}{ll}
	\sum_{j_1=i_1\vee 0}^{n}\sum_{j_2=i_2\vee 0}^{n}\theta^{j_1-i_1}\rho^{j_2-i_2}, \ &\hbox{if } i_1\leq n,\, i_2\leq n,\\
	0, &\hbox{otherwise.}
\end{array}\right.\\
&=\left\{
\begin{array}{ll}
	\frac{ \theta^{ (i_1\vee 0) - i_1} (1 - \theta^{n+1})}{1-\theta} \frac{ \rho^{ (i_2\vee 0) - i_2} (1 - \rho^{n+1})}{1-\rho} , 
	\ &\hbox{if } i_1\leq n,\, i_2\leq n,\\
	0, &\hbox{otherwise}.
\end{array}\right.
\end{align*}
In the above, $i_1 \vee 0 = \max\{i_1, 0\}$. Obviously, the coefficients $\{b_{n i}\}$ are uniformly bounded in $n$ 
and $i = (i_1,i_2)$. To determine the size of $B_n$, we assume $L(\cdot)\equiv1$ for simplicity. Then
\[
\begin{split}
&B_n=\bigg(\sum_{i\in\mathbb{Z}^2}|b_{n i}|^\alpha\bigg)^{1/\alpha} \\
&=\frac 1 { (1-\theta)(1-\rho)} \Bigg(\sum_{i_1= -\infty}^n \sum_{i_2= -\infty}^n \big[ \theta^{ (i_1\vee 0) - i_1}(1-\theta^{n+1})\big]^\alpha 
\big[ \rho^{ (i_2\vee 0) - i_2}(1-\rho^{n+1})\big]^\alpha\Bigg)^{1/\alpha}. \\
\end{split}
\]
We break the summation $\sum_{i_1= -\infty}^n$ into $ \sum_{i_1= -\infty}^{-1} + \sum_{i_1= 0}^n$. Notice that
\[
\sum_{i_1= -\infty}^{-1}   \big[ \theta^{ (i_1\vee 0) - i_1}(1-\theta^{n+1})\big]^\alpha = \frac{\theta^\alpha} 
{1 - \theta^\alpha} (1 - \theta^{n+1})^\alpha
\]
and the second sum is
\[
\sum_{i_1= 0}^n  (1-\theta^{n+1})^\alpha = (n+1)(1-\theta^{n+1})^\alpha.
\]
Hence, we derive  $B_n \sim (1-\theta)^{-1}(1-\rho)^{-1} n^{2/\alpha}$  as $n\to\infty$. 
}\end{example}


\begin{example} {\rm (}{\it Spatial fractional ARIMA models}{\rm)}
{\rm For $d = 2$ and constants  $\beta_1, \beta_2 \in (0, 1)$, as a special case of the following spatial 
fractional ARIMA model:
\begin{equation}\label{Eq:Frima1}
\Phi( \mathbf B_1, \mathbf B_2) X_{i} = \Psi(\mathbf B_1, \mathbf B_2) (1- \mathbf{B}_1)^{-\beta_1} 
(1- \mathbf{B}_2)^{-\beta_2}
\xi_{i},\; i=(i_1, i_2)\in\Z^2, 
\end{equation}
where $\Phi$ and $\Psi$ are polynomials of two variables,  we take $\Phi  \equiv 1$ and $\Psi \equiv 1$. Then 
\begin{equation}\label{Eq:Frima2}
\begin{split}
X_i&=X_{i_1, i_2}= \bigg(\sum_{j_1=0}^\infty\frac{\Gamma(j_1+\beta_1)}{\Gamma(\beta_1) j_1!}\bigg)
\bigg(\sum_{j_2=0}^\infty\frac{\Gamma(j_2+\beta_2)}{\Gamma(\beta_2) j_2!}\bigg) \xi_{i_1,i_2}\\
&= \sum_{(j_1, j_2) \in \mathbb Z_+^2}  a_{j_1, j_2}\xi_{i_1-j_1,i_2-j_2}, 
\end{split}
\end{equation}
where $\Gamma(\cdot)$ is the Gamma function, $\mathbb Z_+^2$ is the set of pairs of non-negative integers 
and the coefficients  
$$
a_{j_1, j_2} = \frac{\Gamma(j_1+\beta_1)
\Gamma(j_1+\beta_2) }{\Gamma(\beta_1) \Gamma(\beta_2)j_1! j_2!}, \quad(j_1, j_2) \in \mathbb Z_+^2.
$$
It follows from Stirling's formula that $a_{j_1, j_2}$  
satisfy
\begin{equation*}\label{Eq:coe-frima}
a_{j_1, j_2} \sim C j_1^{\beta_1-1} j_2^{\beta_2 - 1} \quad \text{as }\, j_1,\, j_2 \to \infty.
\end{equation*}

The linear random field $\{X_{i}, i\in\Z^2 \}$ in \eqref{Eq:Frima1} is an extension of of the fractional ARIMA time 
series with i.i.d. stable innovations considered by Kokoszka and Taqqu (1995) and Kokoszka (1996). A similar 
time series with finite second moment was considered by Dedecker, Merlev\`ede and Peligrad (2011), among 
others. One of the important properties of fractional ARIMA time series is  long-range dependence; see, e.g., 
Beran, et al. (2013) for more information. However, studies of fractional ARIMA random fields
have not been well developed.

In order for \eqref{Eq:Frima2} to be well defined, we assume $\beta_1, \beta_2 \in (0, 1)$ satisfy 
\begin{equation*}\label{Eq:beta}
 (1 - \beta_k)\alpha > 1, \hbox{ for } \, k = 1, 2,
 \end{equation*}
which also requires $\alpha \in (1, 2]$.  

For any integer $n \ge 1$,  
we consider the index set $\Gamma_{n}^2 =  [0, n]^2\cap\Z^2$. 
Then, under the conditions stated above, for every $i \in \mathbb Z^2$,
\[
\begin{split}
b_{ni} &= \sum_{(j_1,j_2) \in \Gamma_{n}^2} a_{j_1-i_1, j_2 - i_2}\\
&\sim \left\{\begin{array}{ll}
C  \sum_{{\tiny 
i_1\vee 0 \le j_1 \le  n, i_2\vee 0 \le j_2 \le n}} ( j_1 - i_1)^{\beta_1-1} (j_2 - i_2)^{\beta_2 - 1},  &\hbox{ if } \, i_1 \le n, i_2 \le  n,\\
0, &\hbox{ otherwise}
\end{array} \right.\\ 
&\sim  \left\{\begin{array}{ll}
\frac{C}{\beta_1\beta_2} \big[( n - i_1)^{\beta_1} - ( - i_1)_+^{\beta_1}\big] \big[(n - i_2)^{\beta_2} - (-i_2)_+^{\beta_2}\big], 
    &\hbox{ if } \, i_1 \le n, i_2 \le  n,\\
0, &\hbox{ otherwise}
\end{array} \right.
\end{split}
\]
as $ n \to \infty.$ In the above, $(-i_1)_+ = 0$ if $i_1 \ge 0$ and $(-i_1)_+ = |i_1|$ if $i_1 < 0$. 

Furthermore, under the additional assumption $L(\cdot)\equiv1$, we can verify that 
\[
\begin{split}
B_n&=\bigg(\sum_{i\in\mathbb{Z}^2}|b_{n i}|^\alpha\bigg)^{1/\alpha}  \\
&  \sim  \frac{C}{\beta_1\beta_2}  \Bigg(\sum_{i_1 \le n }\sum_{i_2 \le n} \big[( n - i_1)^{\beta_1} - ( - i_1)_+^{\beta_1}\big]^\alpha 
\big[(n - i_2)^{\beta_2} - (-i_2)_+^{\beta_2}\big]^\alpha \Bigg)^{1/\alpha} \\
& \sim {C_0}\,n^{\beta_1 + \beta_2 +\frac2 \alpha}, 
\end{split}
\]
as $n \to \infty$, where the constant $C_0$ is given by
\[
C_0 =  \frac{C}{\beta_1\beta_2} \Bigg(\int_{-\infty}^1 \big[(1  - x)^{\beta_1} - (-x)_+^{\beta_1} \big]^\alpha dx\cdot
\int_{-\infty}^1 \big[(1 - x)^{\beta_2} - (-x)_+^{\beta_2}\big]^\alpha dx\Bigg)^{1/\alpha},
\]
where $x_+ = \max\{x, 0\}$ for all $x \in \mathbb R$.
}
\end{example}

The models in the first two examples satisfy the conditions for the main results in Section 3 and 4 
and then we have the convergence in distribution and local limit theorem. 
It is also clear that the two models can be readily extended to the case of $d > 2$.

\begin{example} {\rm (}{\it Linear random fields with isotropic coefficients}{\rm)}
{\rm For any $i\in\mathbb{Z}^d$, let $a_i=\|i\|^{-\beta}$ where $\beta>0$ is a constant and $\|\cdot\|$ denotes the 
Euclidean norm. Obviously, we have that
\begin{itemize}
\item[(i)] $\beta\leq d\Longleftrightarrow \sum_{i\in\mathbb{Z}^d}|a_i|=\infty$,
\item[(ii)] $\alpha\beta>d\Longrightarrow\sum_{i\in\mathbb{Z}^d}|a_i|^\alpha L(1/|a_i|)<\infty$,
\end{itemize}
and $b_{ni}=\sum_{j\in\Gamma_{n}^d}a_{j-i}=\sum_{j\in\Gamma_{n}^d}\|j-i \|^{-\beta}$. 

In order to determine $B_n$,  we assume $L(\cdot)\equiv1$ (for simplicity) and 
$$\Gamma_n^d= \big\{k=(k_1,\ldots,k_d)\in\mathbb{Z}^d:\, |k_l|\leq c_l n, \, l=1,2,\ldots,d \big\},$$
where $c_l > 0$ ($l = 1, \ldots, d$) are constants.  Denote by $I=\prod_{l=1}^d[-c_l,c_l]$ the corresponding rectangle. 
Then we have
\begin{equation}\label{Eq:Bn1}
\begin{split}
B_n&=\bigg(\sum_{i\in\mathbb{Z}^d}|b_{ni}|^\alpha\bigg)^{1/\alpha}
=\bigg(\sum_{i\in\mathbb{Z}^d}\Big|\sum_{j\in\Gamma_{n}^d}\|j-i\|^{-\beta}\Big|^\alpha\bigg)^{1/\alpha}\\
&=n^{d-\beta}\cdot\Bigg(\sum_{i\in\mathbb{Z}^d}\bigg|\sum_{j/n\in I}\|j/n-i/n\|^{-\beta}\cdot n^{-d}\bigg|^\alpha\cdot 
n^{-d}\cdot n^{d}\Bigg)^{1/\alpha}\\
&\sim n^{d-\beta+\frac{d}{\alpha}}\cdot\Bigg(\int_{\mathbb{R}^d}\bigg|\int_I\|x-y\|^{-\beta} dx\bigg|^\alpha dy\Bigg)^{1/\alpha}.
\end{split}
\end{equation}
By splitting the integral in $y$ over $J =  \prod_{l=1}^d[-2c_l,\, 2c_l]$ and $\R^d \backslash J$, we can verify that
\[
\int_{\mathbb{R}^d}\bigg|\int_I\|x-y\|^{-\beta} dx\bigg|^\alpha dy<\infty \, \Longleftrightarrow\, \beta<d\  \text{ and } \ 
\alpha\beta>d.
\]
Therefore, when $\beta<d < \alpha\beta$, we have $B_n\sim C_1n^{(1+\alpha^{-1})d-\beta}$ which shows that the 
effect of ``long memory'' index $\beta$ on the normalizing constants $B_n$.  It follows from the first two lines in 
(\ref{Eq:Bn1}) that  $\sup\limits_i |b_{ni}|\sim C_2 n^{d-\beta}$.
}\end{example}

\begin{example} {\rm (}\it Linear random fields with anisotropic coefficients{\rm)}
{\rm  For any $i\in\mathbb{Z}^d$, let $a_i= \big(\sum_{l=1}^d|i_l|^{\beta_l}\big)^{-\gamma},$ 
where $\beta_1,\ldots,\beta_d$ and $\gamma$ are 
positive constants.  An important feature of this model is that, if $\|i\| \to \infty$ along the $l$th direction of 
 $\mathbb{Z}^d$, then the rate of $a_i \to 0$  is $|i_l|^{-\beta_l \gamma}$.
This class of linear random fields is an extension of the model in Damarackas and  Paulauskas (2017), 
where linear random fields with symmetric $\alpha$-stable innovations were considered.

It is elementary to verify the following statements: 
\begin{equation*}\label{Eq:ani}
\begin{split}
 \sum_{i\in\mathbb{Z}^d}|a_i|=\infty \Longleftrightarrow\int_{[1,\infty)^d}\bigg(\sum_{l=1}^d|x_l|^{\beta_l}\bigg)^{-\gamma}dx_1\ldots dx_d=\infty\,
    \, \Longleftrightarrow\gamma\leq\sum_{l=1}^d\beta_l^{-1}
\end{split}
\end{equation*}
and
\begin{equation*}\label{Eq:ani2}
\alpha\gamma>\sum_{l=1}^d\beta_l^{-1} \Longrightarrow\sum_{i\in\mathbb{Z}^d}\bigg(\sum_{l=1}^d|i_l|^{\beta_l}\bigg)^{-\alpha\gamma} 
L\bigg(\Big(\sum_{l=1}^d|i_l|^{\beta_l}\Big)^{\gamma}\bigg)<\infty.
\end{equation*}

In the following, we consider a family of index sets  $\Gamma_n^d$ which are formed by rectangular lattice points that 
matches the anisotropy in $\{a_i, i \in \mathbb{Z}^d\}$ and determine the rates of $B_n \to \infty$. For simplicity we assume 
again that $L(\cdot)\equiv1$.
 
For any integer $n  \ge 1$, 
 let  $\Gamma_n^d \subset \mathbb Z^d$  be the set of  rectangular lattice points defined by 
\begin{equation*}\label{Eq:G1}
\Gamma_n^d = \Big\{ k \in \mathbb Z^d: |k_l| \le n^{1/\beta_l} \ \hbox{ for } l = 1, \ldots, d\Big\}.
\end{equation*}
Notice that $|\Gamma_n^d| \asymp \prod_{l=1}^d n^{1/\beta_l} = n^Q$, where $Q=\sum_{l=1}^d\beta_l^{-1}$. 
Moreover, the volume of a rectangle of side lengths $n^{-1/\beta_l}$ ($l = 1, \ldots, d$) is $n^{-Q}$.

Now we determine the rate of $B_n \to \infty$. 
Since  
$$b_{ni}=\sum_{j\in\Gamma_{n}^d}a_{j-i} = \sum_{j\in\Gamma_{n}^d}\bigg(\sum_{l=1}^d|j_l-i_l|^{\beta_l}\bigg)^{-\gamma},$$ 
we have
\begin{equation}\label{Eq:G1b}
\begin{split}
B_n&
=\left(\sum_{i\in\mathbb{Z}^d}\bigg|\sum_{j\in\Gamma_{n}^d}\bigg(\sum_{l=1}^d|j_l-i_l|^{\beta_l}\bigg)^{-\gamma}\bigg|^\alpha\right)^{1/\alpha}\\
&=\left(\sum_{i\in\mathbb{Z}^d}\bigg|\sum_{j\in\Gamma_{n}^d}\bigg(\sum_{l=1}^d\bigg|\frac{j_l}{n^{1/\beta_l}}-\frac{i_l}{n^{1/\beta_l}}\bigg|^{\beta_l}\bigg)^{-\gamma}
\cdot n^{-\gamma}\cdot n^{-Q}\cdot n^{Q}\Big|^\alpha\right)^{1/\alpha}\\
&=n^{Q-\gamma+\frac{Q}{\alpha}}\cdot\left(\sum_{i\in\mathbb{Z}^d}\bigg|\sum_{j\in\Gamma_n^d}\bigg(\sum_{l=1}^d\bigg|\frac{j_l}{n^{1/\beta_l}}-\frac{i_l}{n^{1/\beta_l}}\bigg|^{\beta_l}\bigg)^{-\gamma}\cdot n^{-Q}\bigg|^\alpha\cdot n^{-Q}\right)^{1/\alpha}\\
&\sim n^{Q-\gamma+\frac{Q}{\alpha}}\left(\int_{\mathbb{R}^d}\bigg|\int_I\bigg(\sum_{l=1}^d|x_l-y_l|^{\beta_l}\bigg)^{-\gamma} dx\bigg|^\alpha dy\right)^{1/\alpha}
\end{split}
\end{equation}
as $n \to \infty$, in the last line, $I = [-1, 1]^d$. For the last integral in (\ref{Eq:G1b}), it  can be verified that
\begin{equation*}\label{Eq:G1c}
\int_{\mathbb{R}^d}\bigg|\int_I\bigg(\sum_{l=1}^d|x_l-y_l|^{\beta_l}\bigg)^{-\gamma} dx\bigg|^\alpha dy<\infty \ \Longleftrightarrow \
\alpha\gamma>\sum_{l=1}^d\beta_l^{-1}>\gamma.
\end{equation*}
Hence, under the condition $\alpha\gamma>\sum_{l=1}^d\beta_l^{-1}>\gamma$, we have $B_n\sim C_3n^{(1+\alpha^{-1})Q-\gamma}$ 
which also shows that the effect of ``long memory'' on the normalizing constants $B_n$. From the first two lines in (\ref{Eq:G1b}), we 
can derive that $\sup\limits_i|b_{ni}|\sim C_4 n^{Q-\gamma}$.} 
\end{example}


\section{Proofs}

\noindent\textbf{Proof of Theorem \ref{series}}.

The proof is same as the proof of Theorem 1 in Shukri (1976). For readers convenience and completeness, we provide the proof below. 
We can select $\theta$ with $|\theta|<\pi\alpha/2$ $(0<\alpha< 1)$, $\theta=0$ $(\alpha=1)$ and $|\theta|<(2-\pi)\alpha/2$ $(1<\alpha\le 2)$ and adjust the constant $c_\alpha>0$, then rewrite the characteristic function in (\ref{cha1}) as 
\begin{align*}\label{cha2}
\varphi_{\xi}(t)=\exp \big\{-c_\alpha|t|^\alpha {L}(1/|t|) \exp(({\rm sgn} t)\iota \theta)\big\}
\end{align*}
for $t$ in the neighborhood of zero. This is because $|\beta|\le 1$ and 
$$\exp(({\rm sgn} t)\iota \theta)=\cos\theta(1+\iota\tan\theta{\rm sgn} t).$$
Then the characteristic function of $S_n$ can be written as 
\begin{align*}
\varphi_{S_{n}}(t)  &  =\prod_{i\in\mathbb{Z}^d}\varphi_{\xi}(tc_{ni})\\
&  =\exp\left\{  -c_\alpha|t|^{\alpha}\sum_{i\in\mathbb{Z}^d}|c_{ni}|^{\alpha}L(1/|tc_{ni}|)\exp(({\rm sgn} t)({\rm sgn} c_{ni})\iota \theta)\right\}\nonumber\\
&  =\exp\left\{  -c_\alpha|t|^{\alpha}\sum_{i\in\mathbb{Z}^d}|c_{ni}|^{\alpha}L(1/|c_{ni}|)\frac{L(1/|tc_{ni}|)}{L(1/|c_{ni}|)}\exp(({\rm sgn} t)({\rm sgn} c_{ni})\iota \theta)\right\}.\nonumber
\end{align*}
By (\ref{Eq:S1}), (\ref{Eq:S2}), Condition ($A_{2}$) and the definition of slowly varying function at infinity, for any $0<\varepsilon <1$, we can find $N\in \mathbb{N}$, such that for all $n>N$, $c(1-\varepsilon )<\sum_{i\in\mathbb{Z}^d, c_{ni}>0}|c_{ni}|^{\alpha}L(1/|c_{ni}|)<c(1+\varepsilon )$,   $(1-c)(1-\varepsilon )<\sum_{i\in\mathbb{Z}^d, c_{ni}<0}|c_{ni}|^{\alpha}L(1/|c_{ni}|)<(1-c)(1+\varepsilon )$,  and for all $n>N$ and all $i$, we have $1-\varepsilon <L(1/|tc_{ni}|)/L(1/|c_{ni}|)<1+\varepsilon$.
Since $0<\varepsilon <1$ is arbitrary, we have 
\begin{align*}
\lim_{n\rightarrow \infty}\varphi_{S_{n}}(t) 
=\exp\left\{  -c_\alpha c|t|^{\alpha}\exp(({\rm sgn} t)\iota \theta)\right\}\exp\left\{  -c_\alpha (1-c)|t|^{\alpha}\exp(-({\rm sgn} t)\iota \theta)\right\},\nonumber
\end{align*}
and then 
\begin{align*}
S_n\Rightarrow c^{1/\alpha}S' - (1-c)^{1/\alpha} S'',
\end{align*}
 where $S'$ and $S''$ are independent $\alpha$-stable 
random variables that have the same distribution as $S$ in (\ref{Eq:stablelimt}) and the characteristic function as in (\ref{cs}).
\rule{0.5em}{0.5em}

\begin{lemma}\label{bn}
The quantity $B_n$ defined in (\ref{defBn}) satisfies (\ref{limit}).  
\end{lemma}
\begin{proof}
By the definition of $B_{n}$, for $0<\varepsilon <1$ 
\begin{equation*}
\sum_{i\in \mathbb{Z}^{d}}(|b_{ni}|/(B_{n}-\varepsilon ))^{\alpha
}L((B_{n}-\varepsilon )/|b_{ni}|)>1,
\end{equation*}%
and because $x^{2-\alpha }L(x)$ is increasing 
\begin{eqnarray*}
&&\sum_{i\in \mathbb{Z}^{d}}(|b_{ni}|/(B_{n}-\varepsilon ))^{\alpha
}L((B_{n}-\varepsilon )/|b_{ni}|) \\
&=&\sum_{i\in \mathbb{Z}^{d}}(|b_{ni}|/(B_{n}-\varepsilon ))^{2}\left[
((B_{n}-\varepsilon )/|b_{ni}|)^{2-\alpha }L((B_{n}-\varepsilon )/|b_{ni}|)%
\right]  \\
&\leq &\sum_{i\in \mathbb{Z}^{d}}(|b_{ni}|/(B_{n}-\varepsilon)
)^{2}(B_{n}/|b_{ni}|)^{2-\alpha }L(B_{n}/|b_{ni}|) \\
&=&(B_{n}-\varepsilon )^{-2}\sum_{i\in \mathbb{Z}%
^{d}}(|b_{ni}|)^{2}(B_{n}/|b_{ni}|)^{2-\alpha }L(B_{n}/|b_{ni}|).
\end{eqnarray*}%
By letting $\varepsilon \rightarrow 0$ we obtain%
\begin{equation*}
\sum_{i\in \mathbb{Z}^{d}}(|b_{ni}|/B_{n})^{\alpha }L(B_{n}/|b_{ni}|)\geq 1.
\end{equation*}%
On the other hand, by the definition of $B_{n}$ there is a sequence $%
x_{m}=x_{m}(n)$, which is decreasing to $B_{n}$ such that%
\begin{equation*}
\sum_{i\in \mathbb{Z}^{d}}(|b_{ni}|/x_{m})^{\alpha }L(x_{m}/|b_{ni}|)\leq 1.
\end{equation*}%
By the Fatou's lemma and the right continuity of $L,$ by passing with $%
m\rightarrow \infty $%
\begin{equation*}
\sum_{i\in \mathbb{Z}^{d}}(|b_{ni}|/B_{n})^{\alpha }L(B_{n}/|b_{ni}|)\leq
\lim_{m\rightarrow \infty }\inf \sum_{i\in \mathbb{Z}^{d}}(|b_{ni}|/x_{m})^{%
\alpha }L(x_{m}/|b_{ni}|)\leq 1.
\end{equation*}
This proves (\ref{limit}).
\end{proof}

In the stable case, we have the following lemma paralleling to Proposition 2 of Mallik and Woodroofe (2011). 
It is required in the proof of Theorem \ref{thm1'} and Theorem \ref{thm1}. 

\begin{lemma}
\label{lemmarho} For any $p>\max\{1,\alpha\}$ if $0<\alpha<2$, and $p=2$ if $\alpha=2$ and $\E\xi
_{0}^{2}=\infty$, take $q$ with
$1/p+1/q=1$. 
For the linear random field in the form of (\ref{rfsn}), we have
\begin{equation}
\rho_{n}\ll2^{d}[(\Delta_n/{B_n})^{1/(dp+1)}+\Delta{_{n}}/{B_{n}}]. \label{rho1}%
\end{equation}
For the linear random field in the form of (\ref{deflin}), if the sum is over
$\Gamma_{n}^{d}$, which has the form (\ref{defgamma}), then
\begin{equation}\label{Eq:Delta}
\Delta_{n}\ll2^{d}J_{n}^{1/q}
\end{equation}
and hence 
\begin{equation}
\rho_{n}\ll2^{d}(J_{n}^{1/q} /B_{n}%
)^{1/(dp+1)}+2^{2d}J_{n}^{1/q}/B_{n}. \label{rho2}%
\end{equation}
\end{lemma}

\begin{proof}
For the ease on notation, we only prove the results for the case $d=2$. In this case, 
by (\ref{Incr}), 
\begin{equation*}
\Delta b_{i,j}=b_{i,j}-b_{i,j-1}-b_{i-1,j}+b_{i-1,j-1}. \label{incr2}%
\end{equation*}
Suppose
that $\sup_{i\in\mathbb{Z}^{2}}|b_{ni}|=\sup_{r,s\in\mathbb{Z}}|b_{r,s}|$
occurs at $|b_{r_{0},s_{0}}|$. Following the line in Mallik and Woodroofe
(2011), we have
\begin{equation}
b_{r_{0}+r,s_{0}+s}-b_{r_{0},s_{0}+s}-b_{r_{0}+r,s_{0}}+b_{r_{0},s_{0}}%
=\sum_{u=r_{0}+1}^{r_{0}+r}\sum_{v=s_{0}+1}^{s_{0}+s}\Delta b_{u,v}
\label{jump}%
\end{equation}
for $r,s\geq1$. For $m\geq1$, let
\begin{equation}
Q_{m}=\sum_{r=1}^{m}\sum_{s=1}^{m}\sum_{u=r_{0}+1}^{r_{0}+r}\sum_{v=s_{0}%
+1}^{s_{0}+s}|\Delta b_{u,v}|. \label{Q}%
\end{equation}
Putting (\ref{jump}) and (\ref{Q}) together, we have
\begin{equation}
m^{2}|b_{r_{0},s_{0}}|\leq\sum_{r=1}^{m}\sum_{s=1}^{m}(|b_{r_{0}+r,s_{0}%
+s}|+|b_{r_{0},s_{0}+s}|+|b_{r_{0}+r,s_{0}}|)+Q_{m}. \label{relm}%
\end{equation}

First we consider the case $\alpha=2$. By Remark \ref{sv} and Cauchy-Schwarz inequality, we have 
\begin{align*}
\sum_{r=1}^{m}\sum_{s=1}^{m}|b_{r_{0}+r,s_{0}+s}|  &  \leq m\left(  \sum_{r=1}^{m}\sum_{s=1}^{m}|b_{r_{0}+r,s_{0}+s}|^{2}\right)
^{1/2}\\
&  \leq m B_{n}\left(  \sum_{r\in\mathbb{Z}}\sum_{s\in
\mathbb{Z}}(|b_{r,s}|/B_{n})^{2}L(B_n/|b_{r,s}|)\right)  ^{1/2}= mB_n.
\end{align*}

In the case $0< \alpha<2$, for $0<\epsilon<A$ to be decided later, we write
\begin{align*}
\sum_{r=1}^{m}\sum_{s=1}^{m}|b_{r_{0}+r,s_{0}+s}|  &  =\sum_{r=1}^{m}%
\sum_{s=1}^{m}|b_{r_{0}+r,s_{0}+s}|I(B_{n}/|b_{r_{0}+r,s_{0}+s}|>A)\\
&  +\sum_{r=1}^{m}\sum_{s=1}^{m}|b_{r_{0}+r,s_{0}+s}|I(B_{n}/|b_{r_{0}+r,s_{0}+s}|<\epsilon)\\
&  +\sum_{r=1}^{m}\sum_{s=1}^{m}|b_{r_{0}+r,s_{0}+s}|I(\epsilon\leq
B_{n}/|b_{r_{0}+r,s_{0}+s}|\leq A)\\
&  =I+II+III.
\end{align*}

By Hölder's inequality, applied with $p>\max\{1,\alpha\}$ we have for $q$ satisfying
$1/p+1/q=1$
\begin{align*}
I  &  \leq m^{2/q}\left(  \sum_{r=1}^{m}\sum_{s=1}^{m}|b_{r_{0}+r,s_{0}+s}|^{p}I(B_{n}/|b_{r_{0}+r,s_{0}+s}|>A)\right)
^{1/p}\\
&  \leq m^{2/q}B_{n}\left(  \sum_{r\in\mathbb{Z}}\sum_{s\in
\mathbb{Z}}(|b_{r,s}|/B_{n})^{p}I(B_n/|b_{r,s}|>A)\right)  ^{1/p}.
\end{align*}
By the properties of slowly varying functions, there is $A>0$ such that
uniformly in $r$, $s$ and $n$
\[
(B_n/|b_{r,s}|)^{\alpha-p}I(B_n/|b_{r,s}|>A)\leq L(B_n/|b_{r,s}|).
\]
By (\ref{limit}),
\[
I\leq m^{2/q}B_{n}\left(  \sum_{r\in\mathbb{Z}}\sum_{s\in\mathbb{Z}}%
(|b_{r,s}|/B_{n})^{\alpha}L(B_n/|b_{r,s}|)\right)^{1/p}= m^{2/q}B_{n}.
\]

To treat $II$, note that by the definition of $L(x)$ in (\ref{defL}), $L(x)$ is a constant for $x<b$. 
Therefore there exists a sufficiently small $\epsilon>0$, for all $x<\epsilon$ we have%
\[
 x^{-\alpha}L(x)>1.
\]
Then, if $B_n/|b_{r,s}|<\epsilon$, we have $(|b_{r,s}|/B_{n})^{\alpha}L(B_n/|b_{r,s}|)> 1$. 
Recall that $B_{n}$ satisfies  (\ref{limit}). 
Clearly for this
selection of $\epsilon$, we have
\[
II=\sum_{r=1}^{m}\sum_{s=1}^{m}|b_{r_{0}+r,s_{0}+s}|I(B_n/|b_{r_{0}%
+r,s_{0}+s}|<\epsilon)=0.
\]

To treat $III$, 
by the definition of $L(x)$ in (\ref{defL}), we choose $A>b$ and have  
\begin{align*}
L(x)\ge A^{\alpha-2}(2-\alpha)/\alpha:=M_A>0
\end{align*}
for $\epsilon\leq x \leq A$.
So%
\[
L(B_n/|b_{r,s}|)>M_{A}\text{ for }\epsilon\leq B_n/|b_{r,s}|\leq A.
\]
We apply again Hölder's inequality,
\begin{gather*}
III\leq m^{2/q}B_{n}\left(  \sum_{r\in\mathbb{Z}}\sum_{s\in\mathbb{Z}%
}(|b_{r,s}|/B_{n})^{p}I(\epsilon\leq B_n/|b_{r,s}|\leq A)\right)  ^{1/p}\\
\leq m^{2/q}B_{n}\left(  \frac{1}{M_A}\sum_{r\in\mathbb{Z}}\sum
_{s\in\mathbb{Z}}(|b_{r,s}|/B_{n})^{p}L(B_n/|b_{r,s}|)I(\epsilon\leq
 B_n/|b_{r,s}|\leq A)\right)  ^{1/p}\\
\leq m^{2/q}B_{n}\left(  \frac{\epsilon^{\alpha-p}}{M_A}\sum_{r\in
\mathbb{Z}}\sum_{s\in\mathbb{Z}}(|b_{r,s}|/B_{n})^{\alpha}L(B_n/|b_{r,s}|)\right)  ^{1/p}\ll m^{2/q}B_{n}.
\end{gather*}

By the above estimates of $I$, $II$ and $III$ we obtain
\[
\sum_{r=1}^{m}\sum_{s=1}^{m}|b_{r_{0}+r,s_{0}+s}|\ll m^{2/q}B_{n}.
\]

Similarly,
\[
\sum_{r=1}^{m}\sum_{s=1}^{m}|b_{r_{0},s_{0}+s}|\ll m^{1+1/q}B_{n},
\]
and%
\[
\sum_{r=1}^{m}\sum_{s=1}^{m}|b_{r_{0}+r,s_{0}}|\ll m^{1+1/q}B_{n}.
\]
So
\[
m^{2}|b_{r_{0},s_{0}}|\ll m^{1+1/q}B_{n}+Q_{m}.
\]
By (\ref{relm}) and the above considerations we
have
\[
\rho_{n}=\frac{|b_{r_{0},s_{0}}|}{B_{n}}\ll m^{\frac{1}{q}-1}+
\frac{Q_{m}}{m^2B_{n}}.
\]
Notice that in (\ref{Q}),
\[
Q_{m}=\sum_{r=r_{0}+1}^{r_{0}+m}\sum_{s=s_{0}+1}^{s_{0}+m}(r_{0}%
+m+1-r)(s_{0}+m+1-s)|\Delta b_{r,s}|\leq m^{4}\sup_{r,s}|\Delta b_{r,s}|.
\]
Then, recalling our notation $\Delta_{n}=\sup_{r,s}|\Delta b_{r,s}|$,  we
obtain%
\begin{align*}
\rho_{n}  &  \ll m^{\frac{1}{q}-1}+\frac{m^2\Delta_{n}}{B_{n}}
= m^{-\frac{1}{p}}+\frac{m^2\Delta_{n}}{B_{n}}.
\end{align*}
Denoting now by $\lceil x\rceil$ the smallest integer that exceeds $x$ and
letting $m=\left\lceil (B{_{n}/}\Delta_{n}) ^{p/(2p+1)}\right\rceil$ and $M=(B{_{n}/}\Delta_{n})^{p/(2p+1)}$, 
we have
\[
\rho_{n}\ll (\Delta_n/{B_n})^{1/(2p+1)}%
+2(M^2+1)\Delta_{n}/{B_{n}}= 3(\Delta_n/{B_n})^{1/(2p+1)}+2\Delta{_{n}}/{B_{n}}.
\]

In the $d$ dimensional case, 
\begin{align*}
\rho_{n}    \ll  m^{-\frac{1}{p}}+\frac{m^d\Delta_{n}}{B_{n}}.
\end{align*}
Let $m=\left\lceil (B{_{n}/}\Delta_{n}) ^{p/(dp+1)}\right\rceil$ and $M=(B{_{n}/}\Delta_{n}) ^{p/(dp+1)}$, 
we have
\begin{align*}
\rho_{n}&\ll (\Delta_n/{B_n})^{1/(dp+1)}%
+2^{d-1}(M^d+1)\Delta_{n}/{B_{n}}\\
&= (2^{d-1}+1)(\Delta_n/{B_n})^{1/(dp+1)}+2^{d-1}\Delta{_{n}}/{B_{n}}\\
&\ll 2^{d}[(\Delta_n/{B_n})^{1/(dp+1)}+\Delta{_{n}}/{B_{n}}].
\end{align*}
This proves \eqref{rho1}.

Now we prove (\ref{Eq:Delta}) for linear random field defined in (\ref{deflin}%
). By the linearity of $\Delta$ and the definition of $\Gamma_{n}^{2}$, we
have
\begin{align*}
\Delta b_{r,s}  &  =\Delta\bigg[\sum_{(u,v)\in\Gamma_{n}^{2}}a_{u-r,v-s}%
\bigg]=\Delta\bigg[\sum_{k=1}^{J_{n}}\sum_{(u,v)\in\Gamma_{n}^{2}%
(k)}a_{u-r,v-s}\bigg]\\
&  =\sum_{k=1}^{J_{n}}\Delta\bigg[\sum_{(u,v)\in\Gamma_{n}^{2}(k)}%
a_{u-r,v-s}\bigg].
\end{align*}
For each $k\in\{1,2,...,J_{n}\}$, we have
\begin{align*}
\Delta\bigg[  &  \sum\limits_{(u,v)\in\Gamma_{n}^{2}(k)}a_{u-r,v-s}\bigg]\\
&  =a_{\underline{n}_{1}(k)-r,\underline{n}_{2}(k)-s}-a_{(\overline{n}%
_{1}(k)+1)-r,\underline{n}_{2}(k)-s}+a_{(\overline{n}_{1}(k)+1)-r,(\overline
{n}_{2}(k)+1)-s}-a_{\underline{n}_{1}(k)-r,(\overline{n}_{2}(k)+1)-s}.
\end{align*}
See Fortune, Peligrad and Sang (2021) for the detailed calculation of the above formula.
Hence, by Hölder's inequality,%
\begin{align*}
|\Delta b_{r,s}|  &  =\left|\sum_{k=1}^{J_{n}}(a_{\underline{n}_{1}%
(k)-r,\underline{n}_{2}(k)-s}-a_{(\overline{n}_{1}(k)+1)-r,\underline{n}%
_{2}(k)-s}+a_{(\overline{n}_{1}(k)+1)-r,(\overline{n}_{2}(k)+1)-s}%
-a_{\underline{n}_{1}(k)-r,(\overline{n}_{2}(k)+1)-s})\right|\\
&  \leq4J_{n}^{1/q}\Vert a\Vert_{p}\text{ for all }r, s, 
\end{align*}
where $\Vert a\Vert_{p}=(\sum_{i\in{\mathbb{Z}}^{d}}|a_{i}|^{p})^{1/p}$. 
Therefore $\Delta_{n}\ll 4J_{n}^{1/q}$.
This proves \eqref{Eq:Delta} for $d=2$.

Similarly, in the $d$ dimensional case, we have 
\[
\Delta_{n}\ll2^{d}J_{n}^{1/q}, 
\]
hence  \eqref{Eq:Delta} holds. Finally, it is clear that   (\ref{rho2}) 
follows from (\ref{rho1})  and \eqref{Eq:Delta}.
This completes the proof.
\end{proof}

\bigskip
Now we prove the main theorems in Section 4.

\noindent\textbf{Proof of Theorem \ref{general}}.

The proof is based on the study of the characteristic function of the sum
$S_{n}$. According to Lemma 4.5 and arguments in Section VI.4 in Hennion and
Hervé (2001) (see also Theorem 10.7 in Breiman (1992) and Section 10.4 there),
it suffices to prove (\ref{conclusion}) for all continuous complex valued
functions $g$ defined on ${\mathbb{R}}$, $|g|\in L^{1}({\mathbb{R}})$ such
that the transformation%

\[
\hat{g}(t)=\int\nolimits_{{\mathbb{R}}}e^{-\iota tx}g(x)dx
\]
has compact support contained in some finite interval $[-M,M]$.

The inversion formula gives:%
\[
g(x)=\frac{1}{2\pi}\int\nolimits_{{\mathbb{R}}}e^{\iota sx}\hat{g}(s)\ ds.
\]
Therefore, employing a change of variables and taking the expected value%

\[
{{\mathbb{E}}\,}[g(S_{n}+u)]=\frac{1}{2\pi B_{n}}\int\hat{g}\left(  \frac
{t}{B_{n}}\right)  \;\varphi_{S_{n}}\left(  \frac{t}{B_{n}}\right)
\;\exp\left(  \frac{\iota tu}{B_{n}}\right)  \ dt,
\]
where we used the notation
\[
\varphi_{X}(v)={{\mathbb{E}}\,}(\exp(\iota vX)).
\]

By the Fourier inversion formula we also have
\[
f_{\mathcal{L}}(-u)=\frac{1}{2\pi}\int\varphi_{\mathcal{L}}(t)\text{ }\exp(\iota tu)\ dt
\]
and then
\[
f_{\mathcal{L}}\left(  -\frac{u}{B_{n}}\right)  =\frac{1}{2\pi}\int\varphi_{\mathcal{L}}(t)\text{
}\exp\left(  \frac{\iota tu}{B_{n}}\right)  \ dt.
\]
Therefore%
\begin{gather*}
2\pi\left[  B_{n}{{\mathbb{E}}\,}[g(S_{n}+u)]-f_{\mathcal{L}}\left( -\frac{u}{B_{n}%
}\right)  \int_{{\mathbb{R}}}g(x)dx\right] \\
=\int\hat{g}\left(  \frac{t}{B_{n}}\right)  \;\varphi_{S_{n}}\left(  \frac
{t}{B_{n}}\right)  \;\exp\left(  \frac{\iota tu}{B_{n}}\right)  \ dt-\int
\varphi_{\mathcal{L}}(t)\exp\left(  \frac{\iota tu}{B_{n}}\right)  \ dt\int_{{\mathbb{R}}%
}g(x)dx.
\end{gather*}
So
\begin{align*}
&2\pi\left\vert B_{n}{{\mathbb{E}}\,}[g(S_{n}+u)]-f_{\mathcal{L}}\left( -\frac{u}{B_{n}%
}\right)  \int_{{\mathbb{R}}}g(x)dx\right\vert \\
&\leq\int\left|\hat{g}\left(  \frac{t}{B_{n}}\right)  \;\varphi_{S_{n}}\left(
\frac{t}{B_{n}}\right)  -\varphi_{\mathcal{L}}(t)\int_{{\mathbb{R}}}g(x)dx\right|dt.
\end{align*}
By adding and subtracting $\hat{g}(\frac{t}{B_{n}})\varphi_{\mathcal{L}}(t)$, using the
triangle inequality and taking into account that $\hat{g}$ vanishes outside
$[-M,M]$ we bound the last term by%
\begin{align*}
&  \int_{|t|\leq MB_{n}}\left|\hat{g}\left(  \frac{t}{B_{n}}\right)
\;\varphi_{S_{n}}\left(  \frac{t}{B_{n}}\right)  -\hat{g}\left(  \frac
{t}{B_{n}}\right)  \varphi_{\mathcal{L}}(t)\right|dt\\
&  +\int\left|\hat{g}\left(  \frac{t}{B_{n}}\right)  \varphi_{\mathcal{L}}(t)-\varphi
_{\mathcal{L}}(t)\int_{{\mathbb{R}}}g(x)dx \right|dt=I_{n}+II_{n}.
\end{align*}%
\begin{align*}
I_{n}  &  \leq||\hat{g}||_{\infty}\left(  \int_{|t|\leq T}|\varphi
_{\frac{S_{n}}{B_{n}}}(t)-\varphi_{\mathcal{L}}(t){\Large |}dt+\int_{T<|t|\leq MB_{n}%
}|\varphi_{\frac{S_{n}}{B_{n}}}(t)-\varphi_{\mathcal{L}}(t)|dt\right) \\
&  \leq||\hat{g}||_{\infty}\left(  \int_{|t|\leq T}|\varphi_{\frac{S_{n}%
}{B_{n}}}(t)-\varphi_{\mathcal{L}}(t){\Large |}dt+\int_{T<|t|\leq MB_{n}}|\varphi
_{\frac{S_{n}}{B_{n}}}(t)|dt+\int_{T<|t|}|\varphi_{\mathcal{L}}(t){\Large |}dt\right)  .
\end{align*}
Because $S_{n}/B_{n}\Rightarrow \mathcal{L},$ and $\varphi_{\mathcal{L}}(t)$ is integrable,
$I_{n}$ converges to $0$ provided condition (\ref{integral}) is satisfied. The
second term
\[
II_{n}\leq\int\left|\hat{g}\left(  \frac{t}{B_{n}}\right)  -\int_{{\mathbb{R}}%
}g(x)dx \right|\cdot |\varphi_{\mathcal{L}}(t)|dt
\]
converges to $0$ because $\hat{g}$ is continuous and bounded, $g$ and
$\varphi_{\mathcal{L}}$ are integrable and $B_{n}\rightarrow\infty.$ We obtain%

\[
\sup_{u}\left|B_{n}{{\mathbb{E}}\,}[g(S_{n}+u)]-f_{\mathcal{L}}\left( -\frac{u}{B_{n}%
}\right)  \int_{{\mathbb{R}}}g(x)dx\right|\rightarrow0,\text{ as }n\rightarrow\infty
\]
provided
\[
\lim_{T\rightarrow\infty}\lim\sup_{n}\int_{T<|t|<MB_{n}}|\varphi_{\frac{S_{n}%
}{B_{n}}}(t)|dt\rightarrow0.
\]
\rule{0.5em}{0.5em}

\noindent\textbf{Proof of Theorem \ref{thm2}}.

Recall that $\varphi_{S_{n}}(\cdot)$ is the characteristic function of $S_{n}$
and $\varphi_{\mathcal{L}}$ is the characteristic function of $\mathcal{L}$. We have
$\varphi_{S_{n}}(t/B_{n})=\prod_{i\in\mathbb{Z}^{d}}\varphi_{\xi}%
(tb_{ni}/B_{n})$. Let $\gamma_{n}=\sup_{i}|b_{ni}|$.

By Theorem \ref{thm1'}, we have the weak convergence result (\ref{weak}) in
Theorem \ref{general}. Notice that if $\{b_{ni}\}$ are uniformly bounded in
$n$ and $i$, then obviously $\rho_{n}=\gamma_{n}/B_{n}\rightarrow0$ as
$n\rightarrow\infty$. In the case where $\limsup_{n}\sup_{i\in\mathbb{Z}^{d}%
}|b_{ni}|=\infty$, we impose that $\rho_{n}\rightarrow0$.

Now we verify condition (\ref{integral}) in Theorem \ref{general}. A part of
the argument is similar to the proof of Theorem 2 in Shukri (1976). In
(\ref{rfsn}), we renumber $\{b_{ni}\}$ into a sequence $\{b_{k,n}^{\prime
}\}_{k=1}^{\infty}$ such that $\gamma_{n}=|b_{1,n}^{\prime}|\geq
|b_{2,n}^{\prime}|\geq\cdots$. Since (\ref{limit}), we can select a number
$k_{n}\in\mathbb{N}$ such that%
\[
k_{n}=\inf\left\{  l:\sum_{k=1}^{l}(|b_{k,n}^{\prime}|/B_{n})^{\alpha
}L(B_n/|b_{k,n}^{\prime}|)\geq1/2\right\}  .
\]
Then
\[
\frac{1}{2}\leq\sum_{k=1}^{k_{n}}(|b_{k,n}^{\prime}|/B_{n})^{\alpha}%
L(B_n/|b_{k,n}^{\prime}|)<\frac{1}{2}+(\gamma_{n}/{B_{n}})^{\alpha/2}.
\]
Since ${\gamma_{n}}/{B_{n}}\rightarrow0$ as $n\rightarrow\infty$, by the
properties of slowly varying functions, there exists a constant $C>0$, such
that for sufficiently large $n$,
\[
\sum_{k=1}^{k_{n}}(|b_{k,n}^{\prime}|/B_{n})^{\alpha}L(B_n/|b_{k,n}^{\prime}|)\leq C(\gamma_{n}/B_{n})^{\alpha}L(B_n/\gamma_{n})k_{n}.
\]
Hence, for $n$ large enough%
\begin{equation}
k_{n}\geq\frac{1}{2C(\gamma_{n}/B_{n})^{\alpha}L(B_n/\gamma_{n})}.
\label{boundkn}%
\end{equation}
By (\ref{cha1}), for $0<\alpha\le 2$, there exists $\epsilon>0$ such that for all
$|x|\leq\epsilon$, we have
\[
|\varphi_{\xi}(x)|=\exp\{-c_\alpha|x|^{\alpha}L(1/|x|)\}.
\]
Since $|tb_{k,n}^{\prime}/B_{n}|\leq\epsilon$ for $|t|\leq\epsilon B_{n}/|b_{k_{n},n}^{\prime}|$,  
$k>k_n$, we obtain from the above equality
\begin{align}
\left\vert \varphi_{S_{n}}(t/B_{n})\right\vert  &  \leq\prod_{k=k_{n}%
+1}^{\infty}|\varphi_{\xi}(tb_{k,n}^{\prime}/B_{n})|\label{exp}\\
&  =\exp\left\{  -c_\alpha|t|^{\alpha}\sum_{k=k_{n}+1}^{\infty}|b_{k,n}^{\prime
}/B_{n}|^{\alpha}L(B_n/|tb_{k,n}^{\prime}|)\right\}.\nonumber
\end{align}
Now, by Lemma 5.1 in Davydov (1974), without restricting the generality we
shall assume that $L$ also satisfies the following condition
\begin{equation}
\sup_{x>0}\frac{L(x)}{L(ux)}\leq u^{-\alpha/2}\text{ for all }0<u\leq1.
\label{eqiv}%
\end{equation}
Therefore, by applying (\ref{eqiv}) with $u=|T/t|$ we obtain
\[
L\left(  \frac{B_{n}}{|tb_{k,n}^{\prime}|}\right)  =L\left(  \frac{T}{|t|}
\frac{B_{n}}{|b_{k,n}^{\prime}|T}\right)  \geq\left\vert \frac{T}{t}\right\vert
^{\alpha/2}L\left( \frac{B_{n}}{|b_{k,n}^{\prime}|T}\right)  \text{ for
}|t|>T.
\]

By combining inequality (\ref{exp}) with the last inequality we obtain for
$T<|t|\leq\epsilon B_{n}/|b_{k_{n},n}^{\prime}|,$ and for $n$ and $T$
sufficiently large
\begin{align*}
|\varphi_{S_{n}}(t/B_{n})|  &  \leq\exp\left\{-c_\alpha|t|^{\alpha}\sum_{k=k_{n}%
+1}^{\infty}|{b_{k,n}^{\prime}}/{B_{n}}|^{\alpha}L(B_n/|t{b_{k,n}^{\prime}}|)\right\}\\
&  \leq\exp\left\{-c_\alpha T^{\alpha/2}|t|^{\alpha/2}\sum_{k=k_{n}+1}^{\infty}%
|{b_{k,n}^{\prime}}/{B_{n}}|^\alpha L(B_n/|b_{k,n}^{\prime}T|)\right\}\\
&  \leq\exp\left\{-c_\alpha|t|^{\alpha/2}\sum_{k=k_{n}+1}^{\infty}|{b_{k,n}^{\prime}%
}/{B_{n}}|^{\alpha}L(B_n/|b_{k,n}^{\prime}|)\right\}.
\end{align*}
In the last step we used the properties of slowly varying functions. Notice
that $B_n/|b_{k,n}^{\prime}|\rightarrow\infty$ uniformly for all $k$ as
$n\rightarrow\infty$.

By (\ref{limit}) and the selection of $k_{n}$, we have 
\[
\sum_{k=k_{n}+1}^{\infty}|b_{k,n}^{\prime}/B_{n}|^{\alpha}L(B_n/|b_{k,n}^{\prime
}|)\geq  \frac{1}{2}-\left(  \frac{\gamma_{n}}{B_{n}}\right)
^{\alpha/2}
\]
and, as a consequence
\[
\left\vert \varphi_{S_{n}}(t/B_{n})\right\vert \leq\exp\left\{  -c_\alpha |t|^{\alpha
/2}\right\}  .
\]
Then we have found $\epsilon>0\,$such that
\begin{align*}
&  \lim_{n\rightarrow\infty}\int_{T<|t|<\epsilon B_{n}/|b_{k_{n},n}^{\prime}%
|}|\varphi_{S_{n}}(t/B_{n})|dt\\
&  \leq\int_{|t|>T}\exp\left\{  -c_\alpha |t|^{\alpha/2}\right\}  dt,
\end{align*}
which goes to $0$ as $T\rightarrow\infty$. If $\epsilon B_{n}/|b_{k_{n}%
,n}^{\prime}|\geq DB_{n}$, then condition (\ref{integral}) holds and the proof
is complete. Otherwise, it remains to analyze the term%
\[
G_{n}=\int_{\epsilon B_{n}/|b_{k_{n},n}^{\prime}|<|t|\leq DB_{n}}%
|\varphi_{S_{n}}(t/B_{n})|dt.
\]
We shall consider two cases. In the case (i) that $|b_{ni}|$ is uniformly
bounded by a constant $b_0>0$ for all $i\in\mathbb{Z}^{d}$ and $n\in\mathbb{N}$, notice
that $\epsilon\leq\epsilon|b_{k,n}^{\prime}|/|{b_{k_{n},n}^{\prime}}%
|\leq{t|b_{k,n}^{\prime}|}/{B_{n}}\leq Db_0$ for $\epsilon B_{n}/|b_{k_{n}%
,n}^{\prime}|<|t|\leq DB_{n}$ and $1\leq k\leq k_{n}$. Since the innovations
have a non-lattice distribution and $\varphi_{\xi}$ is continuous, there
exists $0<r<1$ such that $|\varphi_{\xi}(tb_{k,n}^{\prime}/B_{n})|\leq r$ for
all $\epsilon B_{n}/|b_{k_{n},n}^{\prime}|<|t|\leq DB_{n}$ and $1\leq k\leq
k_{n}$. Therefore
\begin{align}
G_{n}  &  \leq\int_{\epsilon B_{n}/|b_{k_{n},n}^{\prime}|<|t|\leq DB_{n}}%
\prod_{k=1}^{k_{n}}\left\vert \varphi_{\xi}(tb_{k,n}^{\prime}/B_{n}%
)\right\vert dt\label{limitkn}\\
&  \leq2DB_{n}r^{k_{n}}.\nonumber
\end{align}
By (\ref{boundkn}),
\[
k_{n}\geq\frac{1}{2C(\gamma_{n}/B_{n})^{\alpha-\delta}}>\frac{B_{n}%
^{\alpha-\delta}}{2Cb_0^{\alpha-\delta}}%
\]
for some $0<\delta<\alpha$ and $n$ is sufficiently large. Then we have
$\lim_{n\rightarrow\infty}B_{n}r^{k_{n}}=0$ and $G_{n}\rightarrow0$ as
$n\rightarrow\infty$.

In the case (ii) where $\limsup_{n}\sup_{i\in\mathbb{Z}^{d}}|b_{ni}|=\infty$
we shall use the Cramér condition. By Lemma 5.2 in Fortune, Peligrad and Sang
(2021) we know that under Cramér condition, for every $\epsilon>0$ there
exists $r>0$ such that
\[
|\varphi_{\xi}(t)|\leq r\text{ for all }|t|>\epsilon.
\]
Since $\epsilon|b_{k,n}^{\prime}|/|{b_{k_{n},n}^{\prime}}|\geq\epsilon$ for
$|t|>\epsilon B_{n}/|b_{k_{n},n}^{\prime}|$ and $1\leq k\leq k_{n}$,
$|\varphi_{\xi}(tb_{k,n}^{\prime}/B_{n})|\leq r$ for all $|t|>\epsilon
B_{n}/|b_{k_{n},n}^{\prime}|$ and $1\leq k\leq k_{n}$. Therefore
(\ref{limitkn}) still holds and we also have
\[
\lim_{n\rightarrow\infty}G_{n}=0.
\]
This completes the proof.
\rule{0.5em}{0.5em}

\noindent\textbf{Proof of Theorem \ref{thm3} and Theorem \ref{thm3'}}.

In the case that $\sum
_{i\in\mathbb{Z}^{d}}|a_{i}|<\infty$, $\{b_{ni}\}$ are uniformly bounded on
$i$ and $n$. Hence Theorem \ref{thm3} holds by the first part of Theorem \ref{thm2}.
The proof of Theorem \ref{thm3'} with $1\leq\alpha\leq2$ and $\sum_{i\in
\mathbb{Z}^{d}}|a_{i}|=\infty$ is similar to the proof of the second case in
Theorem \ref{thm2}. Notice that condition $J_{n}=o(B_{n}^{q})$ for $q\geq2$
implies the weak convergence in Theorem \ref{thm1}.
\rule{0.5em}{0.5em}
\vskip10pt


\noindent\textbf{Acknowledgement.} 
The authors are grateful to the referee and the Associate
Editor for carefully reading the paper and for insightful suggestions that significantly
improved the presentation of the paper.
The research of Magda Peligrad is partially supported by NSF grant
DMS-2054598, USA. The research of Hailin Sang is partially supported by the Simons
Foundation Grant 586789, USA. The research of Yimin Xiao is partially supported by 
NSF grant DMS-1855185, USA. Guangyu Yang's work is partially supported by the 
Foundation of Young Scholar of the Educational Department of Henan Province grant 
2019GGJS012, China.
\\


\end{document}